\documentclass[a4wide]{article}
\usepackage[a4paper, total={6in, 8in}]{geometry}
\usepackage[utf8]{inputenc}
\usepackage{amsmath}
\usepackage{amssymb}
\usepackage{mathtools}
\usepackage{mathrsfs}
\usepackage{hyperref}
\usepackage{graphicx}
\usepackage{longtable}
\usepackage{xcolor}
\usepackage{tikz}
\definecolor{axeZ}{RGB}{20,81,204}
\definecolor{axeY}{RGB}{212,43,11}
\definecolor{axeX}{RGB}{136,189,30}
\newcommand\perm[2]{\ensuremath{({\color{axeZ}#1}, {\color{axeY}#2})}}
\newcommand\bsig{{\ensuremath{\boldsymbol{\sigma}}}}
\newcommand\bi{{\ensuremath{\boldsymbol{i}}}}
\newcommand\bj{{\ensuremath{\boldsymbol{j}}}}
\newcommand\bbasig{{\ensuremath{\overline{\boldsymbol{\sigma}}}}}
\newcommand\basig{{\ensuremath{\overline{\sigma}}}}
\newcommand\bapi{{\ensuremath{\overline{\pi}}}}
\newcommand\bpi{{\boldsymbol{\pi}}}
\newcommand\btau{{\boldsymbol{\tau}}}
\newcommand\bdir{{\boldsymbol{dir}}}
\newcommand\Sym{\ensuremath{Sym}}
\newcommand\rev{\text{rev}}
\newcommand\proj{\text{proj}}
\newcommand\Id{Id}
\newcommand\matSym[4]{\ensuremath{\big(\begin{smallmatrix}
	#1 & #2\\
	#3 & #4
\end{smallmatrix}\big)}}
\newcommand\matSymT[9]{\ensuremath{\big(\begin{smallmatrix}
			#1 & #2 & #3\\
			#4 & #5 & #6 \\
			#7 & #8 & #9
\end{smallmatrix}\big)}}

\newenvironment{centermath}
{\begin{center}$\displaystyle}
	{$\end{center}}

\newcommand\vinpat[2]{\ensuremath{#1|_{\color{axeX}#2}}}
\newcommand\vinpatb[3]{\ensuremath{#1|_{{\color{axeX}#2},{\color{axeY}#3}}}}
\newcommand\vinpatd[5]{\ensuremath{\vinpat{\perm{#1}{#2}}{{{\color{axeX}#3},{\color{axeY}#4},
				{\color{axeZ}#5}}}}}

\newcommand\baxpa{\ensuremath{\vinpatb{2413}{2}{2}}}
\newcommand\baxpb{\ensuremath{\vinpatd{312}{213}{1}{2}{.}}}
\newcommand\baxpc{\ensuremath{\vinpatd{3412}{1432}{2}{2}{.}}}
\newcommand\baxpd{\ensuremath{\vinpatd{2143}{1423}{2}{2}{.}}}
\newcommand\sbaxpa{\ensuremath{\Sym(\baxpa)}}
\newcommand\sbaxpb{\ensuremath{\Sym(\baxpb)}}
\newcommand\sbaxpc{\ensuremath{\Sym(\baxpc)}}
\newcommand\sbaxpd{\ensuremath{\Sym(\baxpd)}}

\newcommand\patt{\Sym(\perm{132}{213})}
\newcommand\TODO{{\color{red}TODO}}

\newcommand\permfirst{\perm{253146}{654321}}
\newcommand\permfirstproj{{\color{axeX}264251}}
\newcommand\permbaxter{\perm{14386527}{47513268}}

\usepackage{amsthm}
\usepackage{authblk}

\newtheorem{theorem}{Theorem}[section]
\newtheorem{corollary}{Corollary}[theorem]

\newtheorem{proposition}[theorem]{Proposition}
\newtheorem{remark}[theorem]{Remark}
\theoremstyle{definition}
\newtheorem{definition}{Definition}[section]
\title{Baxter $d$-permutations and other pattern avoiding classes}

\author[1]{Nicolas Bonichon}
\affil[1]{Univ. Bordeaux, CNRS, Bordeaux INP, LaBRI, UMR 5800, F-33400 Talence, France}

\author[1]{Pierre-Jean Morel}

\date{\today}
\begin{document}
\maketitle

\begin{abstract}
	A permutation of size $n$ can be identified to its diagram in which there
	is exactly one point per row and column in the grid $[n]^2$.
	In this paper we consider multidimensional permutations (or
	$d$-permutations), which are identified to their diagrams on the grid
	$[n]^d$ in which there is exactly one point per hyperplane
	$x_i=j$ for $i\in[d]$ and $j\in[n]$.  We first investigate exhaustively all
	small pattern avoiding
	classes. We provide some bijection to enumerate some of these classes and
	we propose some conjectures for others. We then give a
	generalization of well-studied Baxter permutations into this
	multidimensional
	setting. In addition, we provide a vincular pattern avoidance
	characterization
	of Baxter $d$-permutations.
\end{abstract}

\noindent\makebox[\linewidth]{\rule{\textwidth}{0.5px}}

\section{Introduction}\label{sec:intro}
A permutation  $\sigma=\sigma(1),\dots,\sigma(n) \in S_n$ is a bijection from
$[n]:=\{1,2,\dots,n\}$ to itself. The (2 dimensional) \emph{diagram} of
$\sigma$ is simply the set of points $P_\sigma:=\{(i,\sigma(i)), 1\leq i \leq
	n\}$. The diagrams of permutations of size $n$ are exactly the point sets
	such that
each row and column of $[n]^2$ contains exactly one point.

\begin{figure}[!htb]
	\center{\resizebox{0.4\textwidth}{!}{\input{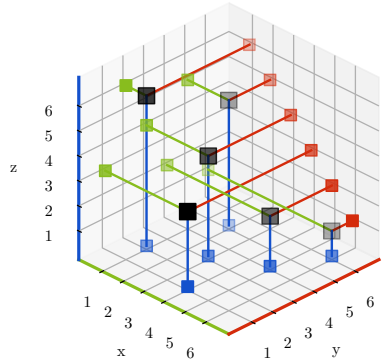}}}
	\caption{The diagram of the 3-permutation \permfirst~together with its
		3 projections
		of dimension 2: the blue, red permutations that define the
		$3$-permutation and green permutation {\color{axeX} $51$} that is
		deduced from the two
		firsts
		permutations.}\label{fig:example1}
\end{figure}

In this paper we are interested in $d$-dimensional diagrams: sets of points
$P_\bsig$ such that every hyperplane $x_i=j$ with $i \in [d]$ and $j \in [n]$
contains exactly one point of $P_\bsig$. Such diagrams are equivalently
described
by a sequence of $d-1$ permutations $\bsig:=(\sigma_1,\dots,\sigma_{d-1})$ such
that
$P_\sigma=\{(i,\sigma_1(i), \sigma_2(i)\dots, \sigma_{d-1}(i)), i \in
	[n]\}$. Figure~\ref{fig:example1} gives an example of a 3-permutation of
size 6. Remark that different generalizations of permutations in higher
dimensions have also been proposed such as Latin
square~\cite{earnest2014permutation, earnest2014permutation} or other
"semi-dense" multidimensional permutations~\cite{eriksson2000combinatorial}.

Permutation-tuple have already been studied (see for
instance~\cite{gunby2019asymptotics,aldred2005permuting}), but as
far as we know, the $d$-permutations have been explicitly considered
only in few papers: \cite{asinowski2010separable,gunby2019asymptotics}.  From
our
point of view the paper of Asinoski and Mansour~\cite{asinowski2010separable}
is the most significant in our context: they present a generalization of
\emph{separable
	permutations} (permutations that can be recursively decomposed with two
elementary composition operations: add the second diagram after the first one
and shift it above or below the first diagram). The formal definition is
provided in
Section~\ref{sec:baxter}. In addition, they characterize these
$d$-permutations with a set of forbidden patterns.

The study of permutations defined by forbidden patterns have received a lot of
attention and sets of small patterns have been exhaustively
studied~\cite{knuth1973art,simion1985restricted,mansour2020enumeration}.
The first main contribution of this paper is to initiate the exhaustive study
of small patterns  for 3-permutations. For this purpose, we propose a
definition of pattern avoidance for $d$-permutations. We say that the
3-permutation $\bsig$ contains the $3$-permutation $\bpi:=\perm{\pi_1}{\pi_2}$
if there is a subset of $P_\bsig$
that is order isomorphic to $P_\bpi$. Also, we say that $\bsig$ contains a
2-permutation $\pi$ if one of its (direct) projections contains $\pi$. We denote
by $S_n^{d-1}(\bpi_1,\dots,\bpi_k)$ the set of $d$-permutations of size $n$
that avoids all patterns $\bpi_1,\dots,\bpi_k$. The
formal definition is provided in Section~\ref{sec:prelim} and is provided for
arbitrary dimensions. This definition is slightly different from the one
introduced in~\cite{asinowski2010separable}. The presented definition has the
advantage to be more expressive than the previous one and it matches with the
classical one for $d=2$.

With this definition in mind, we first investigate exhaustively the enumeration
of 3-permutations defined by small set of patterns to avoid.
Since $3$-permutations are defined by a couple of permutations, it's not
surprising that we fall back on existing combinatorial objects from different
fields: $S_n^2(\perm{12}{12})$ are in bijection with intervals in
the weak-Bruhat order (see Prop~\ref{prop:bruhat}), $S_n^2(\perm{12}{21},
	\perm{312}{132})$ are the
allowable pairs sorted by a priority queue~\cite{atkinson1995priority}. Also,
several "OEIS coincidence" lead us to conjecture other bijections. This is the
case for 4 different pairs of size 3 permutations (see Table~\ref{tab:pat2}).
In
addition, even very simple patterns lead to unknown sequences in the
On-Line Encyclopedia of Integer Sequences, OEIS~\cite{oeis}. This is in
particular the case for all non-trivially equivalent patterns of size 3 (
$S_n^2(\perm{123}{123}), S_n^2(\perm{123}{132}), S_n^2(\perm{132}{213}),
	S_n^2(123), S_n^2(312)$ and $S_n^2(321)$) and some 2 and 3-dimensional
pairs of patterns ($S_n^2(132,\perm{12}{21})$,
$S_n^2(213,\perm{12}{12})$,$S_n^2(231,\perm{12}{12})$,
$S_n^2(231,\perm{21}{12})$,$S_n^2(321,\perm{21}{12})$).

The second main contribution of the paper is a generalization of Baxter
permutations in higher dimension. Baxter permutations are a central family of
permutations that received a lot of attention, in particular because they are
in bijection with a large variety of combinatorial objects: twin binary
trees~\cite{dulucq1998baxter}, plane bipolar
orientations~\cite{bonichon2010baxter}, triples of non-intersecting lattice
paths~\cite{dulucq1998baxter}, Monotone 2-line meanders~\cite{fusy_meanders},
open diagrams~\cite{burrill2016tableau}, Baxter tree-like
tableaux~\cite{aval2021baxter} and boxed arrangements of axis-parallel
segments in $\mathbb{R}^2$~\cite{felsner2011bijections}, and many others.

Having the bijection with boxed arrangements in mind the following
question~\cite{cardinal_open,asinowski2010separable,silveira2018note} was
raised: What is the 3-dimensional analogue of Baxter
permutations? In this paper we propose an analogue of Baxter permutation of
any dimension $d\geq 3$. The
proposed extension seems natural to us, but we didn't investigate the
potential links with boxed arrangements. The generalization of the
bijection with boxed arrangements in higher dimensions remains open.
In addition, we propose a generalization of vincular
patterns for $d$-permutations and we characterize Baxter $d$-permutations
by a set of forbidden vincular patterns (Theorem~\ref{thm:cara-sep}).

The rest of the paper is organized as follows. In Section~\ref{sec:prelim} we
give some definitions and examples of $d$-permutations. We also formalized the
notion of patterns for $d$-permutations and we give few simple properties. Then
in Section~\ref{sec:small-pat} we provide an exhaustive study of the
enumeration of 3-permutations that avoid different sets of small patterns. For
some known sequences, we provide (simple)
explanations. Then in Section~\ref{sec:baxter} we propose a definition of
Baxter $d$-permutations that generalize the classic Baxter permutations. We
also generalize of vincular patterns and we characterize Baxter
$d$-permutations in terms of vincular pattern avoidance. Finally, in
Section~\ref{sec:conclusion} we conclude by a list of open problems.

\section{Preliminaries}\label{sec:prelim}

Let $S_n$ be the symmetric group on $[n] :=\{1,2,\dots,n\}$. Given a permutation
$\sigma=\sigma(1),\dots,\sigma(n) \in S_n$, the \emph{diagram} of $\sigma$,
denoted $P_\sigma$, is the point set $\{(1,\sigma(1)), (2,\sigma(2)), \cdots,
	(n,\sigma(n))\}$. A permutation $\sigma$ \emph{contains} a permutation (or a
\emph{pattern}) $\pi=\pi(1),\dots,\pi(k) \in S_k$ if there exist indices $c_1 <
	\cdots < c_k$ such that $\sigma(c_1) \cdots \sigma(c_k)$ is order-isomorphic to
$\pi$. We say that the set of indices $c_1, \cdots, c_k$ and by extension the
point
set $\{(c_1,
	\sigma(c_1)),\cdots, (c_k, \sigma(c_k))\}$ is an \emph{occurrence}
of the $\pi$.

We denote by $\Id_n$ the identity permutation of size $n$. Given a set of
patterns $\pi_1, \dots, \pi_k$, we denote by $S_n(\pi_1, \dots, \pi_k)$ the set
of permutations of $S_n$ that avoids each pattern $\pi_i$.

\begin{definition}
	A \emph{$d$-permutation} of size $n$, $\bsig:=(\sigma_1,\dots,\sigma_{d-1})$
	is a sequence of $d-1$ permutations of size $n$. We denote by $S_n^{d-1}$
	the set of $d$-permutations of size $n$. Let
	$\bbasig=(\Id_n,\sigma_1,\dots,\sigma_{d-1})$. $d$ is called the
	\emph{dimension} of the permutation.
	The \emph{diagram} of a $d$-permutation $\bsig$ is the set of points in
	$P_\bsig:=\{ (\basig_1(i)$,
	$\basig_2(i),\dots,\basig_{d}(i)), i \in [n] \}$.
\end{definition}

The $d$-permutations of size $n$ are exactly the point sets such that every
hyperplane $x_i=j$ with $i \in [d]$ and $j \in [n]$ contains exactly one point.
One can observe that $|S_n^{d-1}| = n!^{d-1}$.
Figure~\ref{fig:example1} gives an example of a $3$-permutation of size 6.

Given $P:=\{p_1,\dots,p_n\}$ a set of points in $\mathbb{R}^d$ such that every
hyperplane $x_j=\alpha$ with $\alpha\in \mathbb{R}$ contains at most one point
of $P$. The \emph{standardization} of $P$ is the point set
$P'=\{p'_1,\dots,p'_n\}$ in $[n]^{d-1}$ such that the relative order with
respect
to each axis is the same. Hence the standardization of a subset of points of a
diagram is the diagram of a (smaller) $d$-permutation (with the same dimension).

In the sequel we often make the confusion between a $d$-permutation and its
diagram, so that a transformation on one can be directly translated into the
other. For instance, removing a point of a permutation, means removing one
point of its diagram and considering the permutation of the standardization of
the sub-diagram.

\begin{figure}[htb]
	\center{\resizebox{0.9\textwidth}{!}{\input{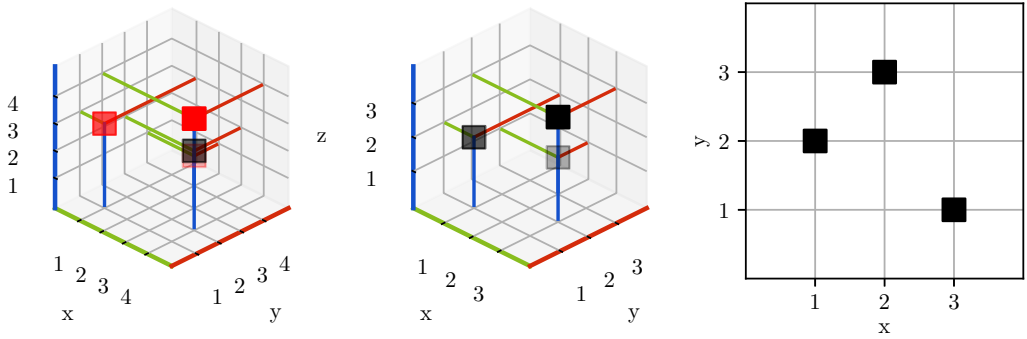}}}
	\caption{On the left, the 3-permutation \perm{1432}{3124}. The red dots are
		an instance of the pattern \perm{132}{213} that is represented in the
		middle. The red dots are also an instance of the pattern 231 that is
		represented on the right.}\label{fig:pattern-ex}
\end{figure}

At this point we are tempted to define a pattern in the following way: a
$d$-permutation $\bsig \in S_n^{d-1}$ contains a pattern $\bpi\in S_k^{d-1}$ if
there exists a subset of points of the diagram of $\bsig$ such that its
standardization is equal to the diagram of $\bpi$ (see
Figure~\ref{fig:pattern-ex}).

This definition have been considered in~\cite{gunby2019asymptotics} for
instance in the context of permutation tuples. For $d=2$, this definition is
consistent with the classical definition over permutations. In higher
dimensions, it is convenient to deal also with patterns of smaller dimensions
(which is not possible when $d=2$). Hence we provide a more general
definition of pattern that matches with the previous one when the dimension of
the pattern is equal to the dimension of the permutation.

Given a sequence of indices $\boldsymbol{i}:=i_1,\dots, i_{d'} \in
[d]^{d'}$, the \emph{projection} on $\boldsymbol{i}$ of $d$-permutation $\bsig$
is the $d'$-permutation
$\proj_{\boldsymbol{i}}(\bsig):=\basig_{i_2}\basig_{i_1}^{-1},
	\basig_{i_3}\basig_{i_1}^{-1},\dots,\basig_{i_{d'}}\basig_{i_1}^{-1}$. $d'$ is
the \emph{dimension} of the projection.

When dealing with permutations of dimension 2 or 3, we often use $x,y,z$ instead of $1,2,3$.

\begin{remark} $\proj_{1,i}(\bsig)=\sigma_{i-1}=\basig_i$,
	$\proj_{i,1}(\bsig)=\basig_i^{-1}$ In particular, when $d=3$,
	$\proj_{xy}(\bsig)={\color{axeZ}\sigma_1}$ and
	$\proj_{xz}(\bsig)={\color{axeY}\sigma_2}$,
	$\proj_{yz}(\bsig)={\color{axeY}\sigma_2}{\color{axeZ}\sigma_1}^{-1}$.
	For instance, $\proj_{yz}(\permfirst)=\permfirstproj$ (see
	Figure~\ref{fig:example1}).
\end{remark}

A projection $\proj_{\boldsymbol{i}}$ is \emph{direct} if $i_1 < i_2 < \dots < i_{d'}$ and \emph{indirect} otherwise.

\begin{definition}
	Let  $\bsig=(\sigma_1, \cdots, \sigma_{d-1}) \in S^{d-1}_n$ -and
	$\bpi=(\pi_1, \cdots, \pi_{d'-1}) \in S^{d'-1}_k$ with $k\leq n$.
	Then $\bsig$ \emph{contains the pattern} $\bpi$,
	if there exist  a direct projection  $\bsig'=\proj_{\boldsymbol{i}}(\bsig)$
	of
	dimension $d'$ and indices $c_1 < \cdots < c_k$ such that
	$\sigma'_i(c_1) \cdots \sigma'_i(c_k)$
	is order-isomorphic to $\pi_i$ for all $i\in [d']$.
	A permutation \emph{avoids} a pattern if it doesn't contain it.
\end{definition}

Given a set of patterns $\bpi_1, \dots, \bpi_k$, we denote by
$S_n^{d-1}(\bpi_1, \dots, \bpi_k)$ the set of $d$-permutations that avoids each
pattern $\bpi_i$.

This definition of pattern differs slightly from the one proposed
in~\cite{asinowski2010separable}: here we consider only \emph{direct}
projections whereas they consider every projections. The advantage of our
convention is that for $d=2$ our definition matches the classical definition of
pattern avoidance: $S_n^2(\bpi) = S_n(\bpi)$, where for instance the set of
$2$-permutations that avoids $2413$ with the other definition is $S_n(2413,
	3142)$, since $3142=\proj_{yx}(2413)$.

We observe that a $d$-permutation $\bsig$ contains a $d$-permutation $\bpi$
if there exists a subset of points of its diagram that have the same relative
positions as those of the diagram of the pattern $\bpi$. This implies that
$\sigma_i \in S(\pi_i) \forall i \in [d-1]$.

Hence
$$S_n(\pi_1)\times S_n(\pi_2) \dots \times S_n(\pi_{d-1}) \subseteq
	S_n^{d-1}(\bpi).$$

In general this inclusion is strict. For instance, the \perm{132}{213} doesn't
contain the pattern \perm{12}{12} but {\color{axeZ} 132} and {\color{axeY} 213}
both contain the pattern 12 (but on different positions).

Avoiding a pattern $\pi$ of dimension 2 means that each projection of dimension $2$ avoids $\pi$, in particular the $d-1$ permutations defining the $d$-permutation, hence

$$S_n^{d-1}(\pi) \subseteq \underbrace{S_n(\pi)\times \dots \times
		S_n(\pi)}_\text{$d-1$ times}.$$

Once again, in general this inclusion is strict. For instance,
$\perm{132}{132} \in S_n(123)\times S_n(123)$ but not in $S_n^2(123)$ since
$\proj_{yz}(\perm{132}{132})={\color{axeX}123}$.

We conclude this section with bijections on $S_n^{d-1}$ that correspond to
symmetries of the $d$-dimensional cube. These operations are defined by
signed permutation matrices of dimension $d$. Let us formalize it. A
\emph{signed permuation matrix} is a square matrix with entries in
$\{-1,0,1\}$ such
that each row and each column contains exactly one non-zero entry. We
denote by \emph{$d$-\Sym} (or simply \emph{\Sym} when the dimension
$d$ is
clear) the set of such matrices of size $d$.

Given $s \in d$-\Sym\ and $\bsig \in S_n^{d-1}$, we define $s(\bsig)$ as
the $d$-permutation whose diagram is the
standardization of the point set $P':=\{(s.(p_1,\dots,p_d)^T)^T,
	(p_1,\dots,p_d) \in P_\bsig \}$.
For instance, in 2 dimension, $\matSym{-1}{0}{0}{1}(\sigma)$ is the
\emph{reverse}
permutation of $\sigma$, denoted by $\rev(\sigma)$:
$\rev(\sigma)(i)=\sigma(n-i+1)$. $\matSym{0}{1}{1}{0}(\sigma)$ is the
\emph{inverse}
permutation of $\sigma$, denoted $\sigma^{-1}$. In dimension 2 there are 8
symmetries and in dimension
3,
there are 48 ($|3$-$\Sym|=48$).

\section{Pattern Avoiding}\label{sec:small-pat}

In this section, we give some exhaustive enumeration of small pattern avoiding
$d$-permutations. We first recall known results for $d=2$ and then we
investigate the case $d=3$.
We start with combinations of basic patterns.
Two sets of patterns $\bpi_1,\bpi_2,...,\bpi_k$ and
$\btau_1,\btau_2,\dots,\btau_{k'}$ are \emph{$d$-Wilf-equivalent} if
$|S_n^{d-1}(\bpi_1,\bpi_2,...,\bpi_k)|=
	|S_n^{d-1}(\btau_1,\btau_2,\dots,\btau_{k'})|$.

We say that 2 sets of
patterns $\bpi_1,\bpi_2,...,\bpi_k$ and
$\btau_1,\btau_2,\dots,\btau_{k'}$ are \emph{trivially $d$-Wilf-equivalent}, if
there
exists a
symmetry $s\in d$-$Sym$ that sends bijectively
$S_n(\bpi_1,\bpi_2,...,\bpi_k)$ on $S_n(\btau_1,\btau_2,\dots,\btau_{k'})$. In
particular, if each pattern $\bpi_1,\bpi_2,...,\bpi_k,
\btau_1,\btau_2,\dots,\btau_{k'}$ is of dimension $d$, the two pattern sets are
equivalent if  $s$ sends the the first one on the second one.

\subsection{Some known result on Permutations}
In dimension 2, there are only 2 patterns of size 2 ($12$ and
$21$) that are trivially Wilf-equivalent. For patterns of size $3$ there are
2 classes of patterns that are trivially Wilf-equivalent: $123$ and $321$
on one hand and $312$, $213$, $231$, $132$ on the other hand. In fact, these
6 patterns are Wilf-equivalent and enumerated by Catalan
numbers~\cite{simion1985restricted}: $|S_n(\tau)|=C_n$ for any $\tau$ of
size $3$ where $C_n=\frac{1}{n+1}\binom{2n}{n}$.
All combinations of patterns of size 3 have been addressed
in~\cite{simion1985restricted} and all combinations size 4
patterns have been studied~\cite{mansour2020enumeration}. Table~\ref{tab:pat2d}
summarizes these results.

\begin{table}[!htb]

	\begin{centermath}
		\begin{array}{|c|c|c|c|}
			\hline
			\text{Patterns}     & \text{\#TWE}                              &
			\text{Sequence}     &
			\text{Comment}                                                                                                              \\
			\hline
			\hline
			12                  & 2                                         & 1, 1, 1, 1, 1, 1, 1, \ensuremath{\cdots}                &
			\\
			\hline
			12, 21              & 1                                         & 1, 0, 0, 0, 0, 0, 0, \ensuremath{\cdots}                &
			\\
			\hline
			312                 & 4                                         & \frac{1}{n+1}\binom{2n}{n}=  1, 2, 5, 14, 42, 132, 429,
			\ensuremath{\cdots} & \text{stack-sortable~\cite{knuth1973art}}                                                             \\
			\hline
			123                 & 2                                         & \frac{1}{n+1}\binom{2n}{n}=  1, 2, 5, 14, 42, 132, 429,
			\ensuremath{\cdots} &
			\cite{knuth1973art}\cite[Prop~19]{simion1985restricted}
			\\
			\hline
			123, 321            & 1                                         & 1, 2, 4, 4, 0, 0, 0, \ensuremath{\cdots}                &
			\cite[Prop~14]{simion1985restricted}
			\\
			\hline
			213, 321            & 4                                         & 1+\frac{n(n-1)}{2} = 1, 2, 4, 7, 11, 16, 22,
			\ensuremath{\cdots} & \cite[Prop~11]{simion1985restricted}
			\\
			\hline
			312, 231            & 2                                         & 2^{n-1} =  1, 2, 4, 8, 16, 32, 64, \ensuremath{\cdots}  &
			\cite[Thm~9]{rotem1981stack}\cite[Prop~8]{simion1985restricted}                                                             \\
			\hline
			231, 132            & 4                                         & 2^{n-1} =  1, 2, 4, 8, 16, 32, 64, \ensuremath{\cdots}  &
			\cite[Prop~9]{simion1985restricted}                                                                                         \\
			\hline
			312, 321            & 4                                         & 2^{n-1} =  1, 2, 4, 8, 16, 32, 64, \ensuremath{\cdots}  &
			\cite[Prop~7]{simion1985restricted}                                                                                         \\
			\hline
			213, 132, 123       & 2                                         & \text{Fibonacci:}~ 1, 2, 3, 5, 8, 13, 21,
			\ensuremath{\cdots} &
			\cite[Prop~15]{simion1985restricted}                                                                                        \\
			\hline
			231, 213, 321       & 8                                         & n = 1, 2, 3, 4, 5, 6, 7, \ensuremath{\cdots}            &
			\cite[Prop~16*]{simion1985restricted}                                                                                       \\
			\hline
			312, 132, 213       & 4                                         & n = 1, 2, 3, 4, 5, 6, 7, \ensuremath{\cdots}            &
			\cite[Prop~16*]{simion1985restricted}                                                                                       \\
			\hline
			312, 321, 123       & 4                                         & 1, 2, 3, 1, 0, 0, 0, \ensuremath{\cdots}                &
			\\
			\hline
			321, 213, 123       & 4                                         & 1, 2, 3, 1, 0, 0, 0, \ensuremath{\cdots}                &
			\\
			\hline
			321, 213, 132       & 2                                         & n = 1, 2, 3, 4, 5, 6, 7, \ensuremath{\cdots}            &
			\cite[Prop~16*]{simion1985restricted}                                                                                       \\
			\hline
		\end{array}
	\end{centermath}
	\caption{Sequences of (2-)permutations avoiding small
		patterns. The second column ($\#TWE$) indicates the number of trivially
		Wilf-equivalent patterns.}\label{tab:pat2d}
\end{table}

\subsection{Exhaustive enumeration of small pattern avoiding 3-permutations}

Here we investigate the different small pattern sets for 3-permutations. We
start in
with combinations of small patterns of dimension 3. The results are synthesized
in Table~\ref{tab:pat1}.

In dimension $3$, they are 4 patterns of size 2 that are trivially
Wilf-equivalent to the pattern \perm{12}{12}. The class $S_n^2(\perm{21}{12})$
corresponds intervals in the weak-Bruhat poset (see
Prop~\ref{prop:bruhat}). An
\emph{inversion} in a
permutation $\pi$ is a pair $(i,j)$
such that $i <j$
and $\pi(i)>\pi(j)$. We say that that a permutation
$\pi_1$ is smaller than a permutation $\pi_2$, $\pi_1 \leq_b \pi_2$, in the
\emph{weak Bruhat order} if the set of inversions of $\pi_1$ is included in
the
set of inversions of $\pi_2$. An \emph{interval} is a pair of comparable
permutations.
No explicit formula is known for the enumeration of intervals in the
weak-Bruhat poset.  This is a
contrast with the 2 dimensional case where almost everything is known for set
of patterns of size at most 4.

\begin{table}[!htb]
	\begin{centermath}
		\begin{array}{|c|c|c|c|}
			\hline
			\text{Patterns}                          & \text{\#TWE} &
			\text{Sequence}                          &
			\text{Comment}
			\\
			\hline
			\hline
			\perm{12}{12}                            &
			4
			                                         &
			1, 3, 17, 151, 1899, 31711, \cdots       &
			\text{Prop~\ref{prop:bruhat}}~
			\href{http://oeis.org/A007767}{A007767}                      \\
			\hline
			\perm{12}{12}, \perm{12}{21}             &
			6
			                                         &
			n! = 1, 2, 6, 24, 120 \cdots             &
			\text{Prop~\ref{prof:factor}}

			\\
			\hline
			\begin{tabular}{@{}c@{}}\perm{12}{12}, \perm{12}{21}, \\
				\perm{21}{12}\end{tabular}                &
			4
			                                         &
			1, 1, 1, 1, 1, 1,
			\cdots                                   &
			\text{Prop~\ref{prof:factor}}

			\\
			\hline
			\begin{tabular}{@{}c@{}}\perm{12}{12}, \perm{12}{21}, \\
				\perm{21}{12}, \perm{21}{21}\end{tabular}                &
			1
			                                         &
			1, 0, 0, 0, 0, 0,\cdots
			                                         &

			\\
			\hline
			\perm{123}{123} & 4 & 1, 4, 35, 524, 11774, 366352, 14953983, \ensuremath{\cdots} & \ensuremath{new} \\
\hline
\perm{123}{132} & 24 & 1, 4, 35, 524, 11768, 365558, 14871439, \ensuremath{\cdots} & \ensuremath{new} \\
\hline
\perm{132}{213} & 8 & 1, 4, 35, 524, 11759, 364372, 14748525, \ensuremath{\cdots} & \ensuremath{new} \\
\hline

			\perm{12}{12}, \perm{132}{312}           &
			48
			                                         &
			(n+1)^{n+1} = 1, 3, 16, 125, 1296 \cdots &
			\href{http://oeis.org/A000272}{A000272}
			\text{\cite{atkinson1995priority,atkinson1993permutational}} \\
			\hline
			\perm{12}{12}, \perm{123}{321}           &
			12
			                                         &
			1, 3, 16, 124, 1262, 15898, \cdots       &
			Prop~\ref{prop:bruhat}~\href{http://oeis.org/A190291}{A190291}
			\\
			\hline
			\perm{12}{12}, \perm{231}{312}           &
			8
			                                         &
			1, 3, 16, 122, 1188, 13844, \cdots       &
			\href{http://oeis.org/A295928}{A295928}? \cite{salo2020cutting}
			\\
			\hline
		\end{array}
	\end{centermath}
	\caption{Sequences of 3-permutations avoiding patterns of dimension 3: 1,2
		or
		3 patterns of size 2 or 1 pattern of size 3.The "?" after
		sequence IDs means that the sequences match on the first terms and that
		we conjecture that the sequences are the same.}\label{tab:pat1}
\end{table}

Avoiding 2 patterns of size 2, also leads to a unique Wilf equivalence class
that is enumerated by $n!$:

\begin{proposition}\label{prof:factor}

	$$|S^2_n(\perm{12}{12}, \perm{12}{12})|=n!,$$
	$$|S^2_n(\perm{12}{12}, \perm{12}{21}, \perm{21}{12})|=1.$$
\end{proposition}
\begin{proof}
	Let us consider the pattern set $\{\perm{12}{21}, \perm{21}{12}\}$ that is
	trivially Wilf equivalent to $\{\perm{12}{12}, \perm{12}{12}\}$.
	Let $\perm{\sigma_1}{\sigma_2}\in S_n^2\{\perm{12}{21}, \perm{21}{12}\}$.
	For all $i,j, \sigma_1(i) < \sigma_1(j)$ if and only if
	$\sigma_1(i) < \sigma_1(j)$. This implies that $\sigma_1=\sigma_2$.
	Hence $S^2_n(\perm{12}{21}, \perm{21}{12})=\{\perm{\sigma}{\sigma}, \sigma
		\in S_n\}$, and  $|S^2_n(\perm{12}{21}, \perm{21}{12})| = n!$.
	In this set, if we avoid a third pattern \perm{21}{21}, the only
	permutation that remains is \perm{\Id_n}{\Id_n}, hence
	$|S^2_n(\perm{12}{21}, \perm{21}{12},\perm{21}{21})|=1$. Since every sets
	of 3 patterns of size 2 are trivially Wilf equivalent, we get the second
	equality.
\end{proof}

As opposed to classical permutations avoiding one pattern of size 3 which are
all
enumerated by Catalan numbers, the patterns of size 3 are not all Wilf
equivalent in dimension 3. Surprisingly, the 3 different classes of
Wilf-equivalent patterns of size 3 lead to new integer sequences.
In contrast, the combination of patterns of size 2 and 3 gives already known
sequences (the link with the last one being only conjectural).

Let us start with the pattern set  $\{\perm{12}{12}, \perm{132}{312}\}$. This
pattern set is sent to the pattern set $\{\perm{12}{21}, \perm{321}{132}\}$ by
the
symmetry \matSymT{0}{0}{-1}{0}{-1}{0}{1}{0}{0}.

The set $S_n^2(\perm{12}{21}, \perm{321}{132})$ are exactly the allowable
pairs sorted by a priority queue as shown in~\cite{atkinson1995priority}.
Moreover it was
proven in~\cite{atkinson1993permutational} that this set is of size
$(n+1)^{n+1}$. A bijection between these permutations and labeled trees has
been described in~\cite{atkinson1995priority}.

\begin{proposition}\label{prop:bruhat}
	\begin{enumerate}
		\item $S_n^2(\perm{12}{12})$ is in bijection with intervals in
		      weak-Bruhat
		      poset.
		\item $S_n^2(\perm{12}{12}, \perm{123}{321})$ is in bijection with
		      intervals in weak-Bruhat that are distributive lattices.
	\end{enumerate}
\end{proposition}
\begin{proof}
	\begin{enumerate}
		\item  Observe that $i_1,i_2$ is an inversion in
		      $\pi_1$
		      but not in $\pi_2$, then $i_1,i_2$ is an instance of the pattern
		      $\perm{12}{12}$ in \perm{\pi_1}{\pi_2}. Hence the class
		      $S_n^2(\perm{21}{12})$
		      corresponds to intervals in the weak-bruhat poset. We conclude by
		      observing that The symmetry $\matSymT{-1}{0}{0}{0}{1}{0}{0}{0}{-1}$
		      sends
		      $S_n^2(\perm{21}{12})$ on $S_n^2(\perm{12}{12})$

		\item As shown in~\cite[Proposition 2.3]{stembridge1996fully}, the
		      sub-poset
		      defined by the interval $\sigma_1,\sigma_2$ is isomorphic to the
		      sub-poset of permutations smaller than $\sigma_1^{-1} \sigma_2$.
		      Moreover, as shown in~\cite[Theorem
			      3.2]{stembridge1996fully}, this sub-poset is a distributive
			      lattice if and only if
		      $\sigma_1^{-1} \sigma_2 \in S_n(321)$. Let $G_n$ be the set of
		      3-permutations $\bsig \in
			      S_n^2(\perm{21}{12})$,
		      such that $\sigma_1^{-1} \sigma_2 \in S_n(321)$.	Let us show that
		      $S_n^2(\perm{21}{12}, \perm{123}{321})= G_n$. If $i_1 < i_2 < i_3$ is
		      an occurrence of $\perm{123}{321}$ in a permutation $\bsig$, then it is
		      also an occurrence of 321 in $\sigma_1^{-1} \sigma_2$. Hence $G_n
			      \subseteq
			      S_n^2(\perm{21}{12}, \perm{123}{321})$ so let us
		      focus on the second inclusion. Let us
		      consider $\perm{\sigma_1}{\sigma_2} \in  S_n^2(\perm{21}{12})$ such
		      that $i_1 < i_2 <i_3$ is an occurrence of 321 in
		      $\sigma_1^{-1} \sigma_2$. if $\sigma_1(i_1) < \sigma_1(i_2)$, then
		      $i_1,i_2$
		      is an occurrence of $\perm{21}{12}$ in $\bsig$, which is
		      impossible. Hence $\sigma_1(i_1) >
			      \sigma_1(i_2)$ Applying the same argument on $i_2$ and $i_3$, we get
		      that $i_1,i_2,i_3$ is an occurrence of 123 in $\sigma_1$.
		      $\sigma_1^{-1} \sigma_2$ and $\sigma_1$ fully determine $\sigma_2$ and
		      we have $\pi_2(i_1) > \pi_2(i_2) > \pi_2(i_3)$. Hence $i_1, i_2, i_3$
		      is an occurrence of $\perm{123}{321}$ in
		      $\bsig$, which concludes the second inclusion.

		      We conclude by observing that the symmetry
		      $\matSymT{-1}{0}{0}{0}{1}{0}{0}{0}{-1}$
		      sends bijectively $S_n^2(\perm{21}{12}, \perm{123}{321})$ on
		      $S_n^2(\perm{12}{12}, \perm{123}{321})$.
	\end{enumerate}
\end{proof}

Now, let us focus on 3-permutations that avoid patterns of dimension 2.
Table~\ref{tab:pat2} synthesizes the results. We start by some considerations on
the trivially $d$-Wilf-equivalence of patterns (and pattern sets) of smaller
dimension.

\begin{remark}\label{rem:triv} Let $\bsig \in S_n^2$ with $n\geq 2$. One can
	observe that if $\proj_{x,y}(\bsig) \in S_n(21)$ and $\proj_{x,z}(\bsig)
		\in
		S_n(21)$ then $\proj_{y,z}(\bsig)$ contains the pattern $21$. Hence
	$|S_n^2(21)|=0$ for $n\geq 2$. On the other hand, one can check that
	$S_n^2(21)=\{\perm{\Id_n}{\Id_n}\}$. More generally, 2 patterns of dimension
	$d$
	can be trivially $d$-Wilf-equivalent but not $d'$-Wilf-equivalent for
	$d'>d$.
	For instance, $12$ and $21$ are trivially 2-Wilf-equivalent but not
	3-Wilf-equivalent. In fact, any symmetry of the 3-cube other than the
	identity
	sends the pattern $12$ into the pattern set $\{12,21\}$.
\end{remark}

Given a symmetry $s\in d$-$Sym$ and an increasing sequence of indices $i_1 <
i_2\dots i_{d'}$,  we define $s_\bi$ as an
element of $d'$-$Sym$ obtained from $s$ by keeping the rows of index in $\bi$
and column containing a no-zero value on one of these rows. For instance if $s
=  \matSymT{0}{0}{-1}{0}{-1}{0}{1}{0}{0}$ and $\bi = 1,3$ then $s_\bi =
\matSym{0}{-1}{1}{0}$.
Given $s \in d$-$Sym$ and $\bpi \in S_n^{d'-1}$, we define as follows, if
$\bpi$
is a $d'$-multipermutation,
$\widetilde{s}(\{\bpi\}):=\{s_\bi(\bpi), \bi =i_1,\dots, i_{d'}\}$ and if
$\bpi_1, \dots, \bpi_k$ is a set, $\widetilde{s}(\{\bpi_1, \dots,
\bpi_k\}):=\cup_{i=1}^k \widetilde{s}(\{\bpi_i\})$.

In general $\widetilde{s}(\widetilde{s^-1}(\bpi)) \neq \bpi$. For instance, as
we saw above, for
$d=3$ and $s$ is the identity matrix of size 3,
$\widetilde{s^{-1}}(\widetilde{s}(\{12\}))= \{12,21\}$.

\begin{proposition}\label{prop:TWE}
	 Two pattern sets $\bpi_1, \dots, \bpi_k$ and $\btau_1,\dots,
	 \btau_k'$ are trivially $d$-Wilf-equivalent if there exists $s\in
	 d$-$Sym$ such that
	 $\widetilde{s}(\bpi_1, \dots, \bpi_k) = \btau_1,\dots, \btau_k'$ and
	 $\bpi_1,
	 \dots, \bpi_k = \widetilde{s^{-1}}(\btau_1,\dots, \btau_k')$.
\end{proposition}
\begin{proof}
	Let $\bpi_1, \dots, \bpi_k$, $\btau_1,\dots, \btau_k'$ and $s$ be as in the
	proposition. Let us first show that $|S_n(\bpi_1, \dots, \bpi_k)| \geq
	|S_n(\btau_1,\dots, \btau_k)|$ and then we will show the other inequality.

	Let $\bsig \not \in S_n^d(\bpi_1, \dots, \bpi_k)$ and let $\bi,k$ such that
	$Proj_\bi(\bsig)$ contains $\bpi_k$. Hence $s_\bi(Proj_\bi(\bsig)))$
	contains $s_\bi(\bpi_k)$.
	Let $\bj$ the indices of the columns the contains a non-zero entry in the
	rows of index in $\bi$ in $s$. Since
	$Proj_\bj(s(\bsig)) = s_\bi(Proj_\bi(\bsig))$ and $s_\bi(\bpi_k) \in
	\widetilde{s}(\bpi_k) \subset
	\{\btau_1,\dots, \btau_k'\}$, we
	have $s(\bsig) \not \in
	S_n^d(\btau_1, \dots, \btau_k)$. Hence $|S_n(\bpi_1, \dots, \bpi_k)| \geq
	|S_n(\btau_1,\dots, \btau_k)|$.

	We proceed similarly for the other inequality: if
	Let $\bsig \not \in S_n^d(\btau_1, \dots, \btau_k)$ and let $\bi,k$ such
	that
	$Proj_\bi(\bsig)$ contains $\btau_k$. Hence $s^{-1}_\bi(Proj_\bi(\bsig)))$
	contains $s^{-1}_\bi(\btau_k)$.
	Let $\bj$ the indices of the columns the contains a non-zero entry in the
	rows of index in $\bi$ in $s^{-1}$. Since
	$Proj_\bj(s^{-1}(\bsig)) = s^{-1}_\bi(Proj_\bi(\bsig))$ and
	$s^{-1}_\bi(\btau_k) \in
	\widetilde{s^{-1}}(\btau_k) \subset
	\{\bpi_1,\dots, \bpi_k'\}$, we
	have $s(\bsig) \not \in
	S_n^d(\bpi_1, \dots, \bpi_k)$. Hence $|S_n(\bpi_1, \dots, \bpi_k)| \leq
	|S_n(\btau_1,\dots, \btau_k)|$.

\end{proof}

What is very surprising is that all classes composed of a single pattern of size
3 lead to a new sequences and that 4 out of 5 classes composed of pairs of
patterns of size 3 seem to match with known sequences. For known sequences
we
didn't find any simple interpretations. If we now consider combination of
patterns of dimension 2 and 3 (see Table~\ref{tab:pat3}),
we find several finite sets, 2 new sequences and 5 sequences that seem to
match with known sequences. 3 out of the 4 couples of patterns of size 2
are in fact equivalent to a single pattern (12 or 21), since any instance
of the pattern of dimension 3 is also an instance of the pattern of
dimension 2.

\begin{table}[!htb]
	\begin{centermath}
		\begin{array}{|c|c|c|c|}
			\hline
			\text{Patterns} & \text{\#TWE} &
			\text{Sequence} &
			\text{Comment}                                                                              \\
			\hline
			\hline
			12              & 1            & 1, 0, 0, 0, 0, \ensuremath{\cdots} & Remark~\ref{rem:triv} \\
			\hline
			21              & 1            & 1, 1, 1, 1, 1, \ensuremath{\cdots} & Remark~\ref{rem:triv} \\
			\hline
			123 & 1 & 1, 4, 20, 100, 410, 1224, 2232, \ensuremath{\cdots} & \ensuremath{new} \\
\hline
132 & 2 & 1, 4, 21, 116, 646, 3596, 19981, \ensuremath{\cdots} & \ensuremath{new} \\
\hline
231 & 2 & 1, 4, 21, 123, 767, 4994, 33584, \ensuremath{\cdots} & \ensuremath{new} \\
\hline
321 & 1 & 1, 4, 21, 128, 850, 5956, 43235, \ensuremath{\cdots} & \ensuremath{new} \\
\hline

			123, 132 & 2 & 1, 4, 8, 8, 0, 0, 0, \ensuremath{\cdots} &  \\
\hline
123, 231 & 2 & 1, 4, 9, 6, 0, 0, 0, \ensuremath{\cdots} &  \\
\hline
123, 321 & 1 & 1, 4, 8, 0, 0, 0, 0, \ensuremath{\cdots} &  \\
\hline
132, 213 & 1 & 1, 4, 12, 28, 58, 114, 220, \ensuremath{\cdots} & \ensuremath{new} \\
\hline
132, 231 & 4 & 1, 4, 12, 32, 80, 192, 448, \ensuremath{\cdots} & \href{http://oeis.org/A001787}{A001787}? \\
\hline
132, 321 & 2 & 1, 4, 12, 27, 51, 86, 134, \ensuremath{\cdots} & \href{http://oeis.org/A047732}{A047732}? \\
\hline
231, 312 & 1 & 1, 4, 10, 28, 76, 208, 568, \ensuremath{\cdots} & \href{http://oeis.org/A026150}{A026150}? \\
\hline
231, 321 & 2 & 1, 4, 12, 36, 108, 324, 972, \ensuremath{\cdots} & \href{http://oeis.org/A003946}{A003946}? \\
\hline

		\end{array}
	\end{centermath}
	\caption{Sequences of 3-permutations avoiding at most 2 patterns of size 2
		or
		3 of dimension 2.The "?" after
		sequence IDs means that the sequences match on the first terms and that
		we conjecture that the sequences are the same.}\label{tab:pat2}
\end{table}

\begin{table}[!htb]
	\begin{centermath}
		\begin{array}{|c|c|c|c|}
			\hline
			\text{Patterns}           & \text{\#TWE} &
			\text{Sequence}           &
			\text{Comment}                             \\
			\hline
			\hline
			12, \perm{12}{12}         &
			1                         & 1, 0, 0,
			0, 0, \ensuremath{\cdots} & 12             \\
			\hline
			12, \perm{21}{12}         &
			3                         & 1, 0, 0,
			0, 0, \ensuremath{\cdots} & 12             \\
			\hline
			21, \perm{12}{12}         &
			1                         & 1, 0, 0,
			0, 0, \ensuremath{\cdots} &                \\
			\hline
			21, \perm{21}{12}         &
			3                         & 1, 1, 1,
			1, 1, \ensuremath{\cdots} & 21             \\
			\hline
			123, \perm{12}{12} & 1 & 1, 3, 14, 70, 288, 822, 1260, \ensuremath{\cdots} & \ensuremath{new} \\
\hline
123, \perm{12}{21} & 3 & 1, 3, 6, 6, 0, 0, 0, \ensuremath{\cdots} &  \\
\hline
132, \perm{12}{12} & 2 & 1, 3, 11, 41, 153, 573, 2157, \ensuremath{\cdots} & \href{http://oeis.org/A281593}{A281593}? \\
\hline
132, \perm{12}{21} & 6 & 1, 3, 11, 43, 173, 707, 2917, \ensuremath{\cdots} & \href{http://oeis.org/A026671}{A026671}? \\
\hline
231, \perm{12}{12} & 2 & 1, 3, 9, 26, 72, 192, 496, \ensuremath{\cdots} & \href{http://oeis.org/A072863}{A072863}? \\
\hline
231, \perm{12}{21} & 4 & 1, 3, 11, 44, 186, 818, 3706, \ensuremath{\cdots} & \ensuremath{new} \\
\hline
231, \perm{21}{12} & 2 & 1, 3, 12, 55, 273, 1428, 7752, \ensuremath{\cdots} & \href{http://oeis.org/A001764}{A001764}? \\
\hline
321, \perm{12}{12} & 1 & 1, 3, 2, 0, 0, 0, 0, \ensuremath{\cdots} &  \\
\hline
321, \perm{12}{21} & 3 & 1, 3, 11, 47, 221, 1113, 5903, \ensuremath{\cdots} & \href{http://oeis.org/A217216}{A217216}? \\
\hline

		\end{array}
	\end{centermath}
	\caption{Sequences of 3-permutations avoiding a permutation of size 2 and
		dimension 3 with a pattern of dimension 2 of size 2 or 3. The "?" after
		sequence IDs means that the sequences match on the first terms and that
		we conjecture that the sequences are the same.}\label{tab:pat3}
\end{table}

We conclude this section with sets of patterns that are invariant by all
symmetries.
Given a $d$-permutation $\bsig$, we denote by $\Sym(\bsig):=\{s(\bsig)) | s \in
	d$-$\Sym\}$.

Figure~\ref{fig:patt3} describes all the symmetric $2$-permutations obtained
from $\perm{132}{213}$. This symmetric pattern plays an important role in
separable $d$-permutations and Baxter $d$-permutations as we will see in
Section~\ref{sec:baxter}.
\begin{remark}\label{rem:patt}
	A convenient way to describe this pattern is the following: a permutation
	$\bsig$ contains the pattern $\patt$ if its diagram contains 3 points
	$p_1,p_2,p_3$ and 3 axes such that $p_1$ and $p_2$ are in the same
	quadrant of $p_3$ in the plane generated by the 2 first axes and $p_3$ is
	between $p_1$ and $p_2$  on the third axis.
\end{remark}
\begin{figure}[htb]
	\center{\input{patt3.pgf}}
	\caption{The 8 3-permutations of \patt.}\label{fig:patt3}
\end{figure}


\begin{table}[!htb]
	\begin{centermath}
		\begin{array}{|c|c|c|c|}
			\hline
			\text{Patterns}              & |\Sym(\bpi)|   &
			\text{Sequence}              & \text{Comment}         \\
			\hline
			\hline
			\Sym(\perm{123}{123})        & 4              & 1, 4,
			32, 368, 4952, 68256, \cdots & new                    \\
			\hline
			\Sym(\perm{123}{132})        & 24             & 1, 4,
			12, 4, 4, 4, \cdots          &
			\text{Prop~\ref{prop:const}}
			\\
			\hline
			\Sym(\perm{132}{213})        & 8              & 1, 4,
			28, 256, 2704, 31192, \cdots & new                    \\
			\hline

		\end{array}
	\end{centermath}
	\caption{Sequences of 3-permutations avoiding a pattern of size 3 with all
		its symmetries. The second column indicates the number of forbidden
		patterns.}\label{tab:pat3b}
\end{table}

The number of permutations avoiding $\Sym(\perm{123}{132})$ become constant (equals to 4) for size greater than 4).
In fact, it can be shown that this permutations are 4 diagonals of the cube.

\begin{proposition}\label{prop:const}
	\[
		S^2_n(\Sym(\perm{123}{132})) =
		\begin{cases}
			S_n^2                                 & \text{if } n \leq 2 \\
			S_3^2 \setminus \Sym(\perm{123}{132}) & \text{if } n = 3    \\
			{\Sym(\perm{Id_n}{Id_n})}             & \text{else}         \\
		\end{cases}
	\]
\end{proposition}

\begin{proof} For $n\leq 4$ the proposition can be easily checked manually.
	For $n \geq 4$ let us show that $S^2_n(\Sym(\perm{123}{132}))=
		\Sym(\perm{Id_n}{Id_n}) =
		\{\perm{Id_n}{\Id_n},\perm{Id_n}{\rev(Id_n)},\perm{\rev(Id_n)}{Id_n},\perm{\rev(Id_n)}{\rev(Id_n)}\}$.
	Clearly, $\Sym(\perm{Id_n}{Id_n}) \subseteq S^2_n(\Sym(\perm{123}{132}))$,
	so let us show the other inclusion.

	Suppose that the proposition is true until  some $n \geq 4$ and let's
	show that it is still true  for $n+1$. Let $\bsig \in
		S^2_{n+1}(\Sym(\perm{123}{132}))$. Let $\bsig'$ the permutation
	obtained by removing the point $(x,y,z)$ such that $z=n+1$. If $\bsig$
	avoids a pattern $\bpi$, $\bsig'$ also avoids $\bpi$. Hence $\bsig' \in
		S^2_{n}(\Sym(\perm{123}{132}))$. By our inductive hypothesis, $\bsig'
		\in \Sym(\perm{Id_n}{Id_n})$. Let us just show that if $\bsig' =
		\perm{Id_n}{Id_n}$, then $\bsig=\perm{Id_{n+1}}{Id_{n+1}}$, the 3 other
	cases being equivalent.
	Let us consider all the different possible position for the point $(x,y,n+1)$. Here we only consider cases where $x \leq y$, the other cases being deduced from the first ones by symmetry:
	\begin{itemize}
		\item $x=y=n+1$. In this case $\bsig=\perm{Id_{n+1}}{Id_{n+1}}$.

		\item $x = y = 1$, the permutation will be $\bsig=\perm{Id_{n+1}}{(n+1) \ 1 \cdots n}$ which contains the pattern $\perm{123}{312} \in \Sym(\perm{123}{132})$. Contradiction.

		\item $x = 1, y > 1$.  $\perm{y \ 1 \cdots y-1 \ y+2 \ \cdots n+1}{n+1 \ 1 \cdots n}$ which contains $\perm{123}{312} \in \Sym(\perm{123}{132})$. Contradiction.

		\item $1 < x < n+1, \ y=x$.
		      $\bsig=\perm{Id_{n+1}}{1 \cdots (x-1) \ (n+1) \ x \cdots n}$ which contains the pattern $\perm{123}{132} \in \Sym(\perm{123}{132})$. Contradiction.

		\item $1 < x < n+1, \ y > x$.
		      $\bsig=\perm{1 \cdots (x-1) \ y \ x \cdots (n+1)}{1 \cdots (y-1) \ (n+1) \ y \cdots n}$  contains $\perm{132}{132} \in \Sym(\perm{123}{132})$. Contradiction.

		\item $x = n+1, \ y < n+1$.
		      $\bsig=\perm{1 \cdots (y-1) (y+1)  \cdots (n+1) \ y}{Id_{n+1}}$ which contains $\perm{231}{123} \in \Sym(\perm{123}{132})$. Contradiction.
	\end{itemize}
	So if $\bsig'=\perm{\Id_n}{\Id_n}$ that $\bsig=\perm{\Id_{n+1}}{\Id_{n+1}}$. By symmetry, we conclude that $\Sym(\perm{Id_{n+1}}{Id_{n+1}}) = S^2_{n+1}(\Sym(\perm{123}{132}))$. Hence the property is true for all $n\geq 4$.
\end{proof}.

We give in appendix sequences corresponding to larger patterns. At that date,
none of these sequences appears in OEIS~\cite{oeis}.

\section{Baxter $d$-permutations}\label{sec:baxter}

In this section we consider Separable $d$-permutations and Baxter $d$-permutations. We
first recall definitions and properties in the classical case ($d=2$). Then we recall the definition and characterization of separable
$d$-permutations given in~\cite{asinowski2010separable}, and after we propose a
definition of Baxter $d$-permutation and we show how some of the properties of
Baxter permutations are generalized in higher dimension. Finally we show that we can also extends the notion of \emph{complete Baxter permutations} and \emph{anti-Baxter permutations}.
\subsection{Separable permutations and Baxter permutations}

Let $\sigma$ and $\pi$ two permutations respectively of size $n$ and $k$. Their
\emph{direct sum} and \emph{skew sum} are the permutations of size $n+k$
defined by:
$$\sigma \oplus \pi := \sigma(1),\dots,\sigma(n),\pi(1)+k,\dots,\pi(k)+n
	\text { and}$$
$$\sigma \ominus \pi := \sigma(1)+k,\dots,\sigma(n)+k,\pi(1),\dots,\pi(k).$$
A permutation is \emph{separable} if it is of size $1$ or it is the direct
sum or the skew sum of two separable permutations. Let us denote by
$Sep_n$ the set of separable permutations of size $n$. These permutations
are enumerated by large Schröder numbers as shown
in~\cite{brightwell1992random}:
$$|Sep_n| =
	\frac1{n-1}\sum_{k=0}^{n-2}\binom{n-1}k\binom{n-1}{k+1}(-1)^{n-k-1}.$$

The characterization of separable permutations with patterns have been
given in~\cite{bose1998pattern}:
$$Sep_n=S_n(2413, 3142).$$

\begin{figure}[htb]
	\center{\resizebox{0.6\textwidth}{!}{\input{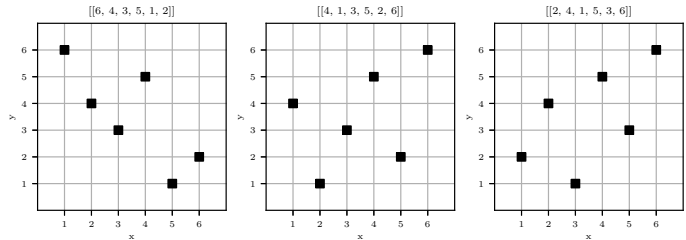}}}
	\caption{On the left the separable permutation $643512 = 1 \ominus((1\ominus1)\oplus1)\ominus(1\oplus1)$. In the middle a Baxter permutation that is not a separable permutation. On the right a permutation that is not a Baxter permutation.}\label{fig:2dsep}
\end{figure}

A related and richer class of permutations are \emph{Baxter permutations}. To
introduce them we first need to define a more general type of patterns.

A vincular pattern is a pattern where some entries must be consecutive in the
permutation.  More formally, a \emph{vincular pattern} $\vinpat{\pi}{X}$ is
composed of $\pi \in S_k$ a permutation and ${\color{axeX}X} \subseteq [k-1]$
a set of (horizontal) \emph {adjacencies}. A permutation $\sigma \in S_n$
\emph{contains} the vincular pattern $\vinpat{\pi}{X}$, if there exist indices
$i_1 < \dots < i_k$ such that $\sigma_{i_1},\sigma_{i_2}\dots\sigma_{i_k}$ is
an
occurrence of $\pi$ in $\sigma$ and that $i_{j+1}=i_{j}+1$ for any $j\in
{\color{axeX}X}$. A vincular pattern $\vinpat{\pi}{X}$ is classically
represented as a permutation with dash between entries without adjacency
constraints. For instance, the vincular pattern $\vinpat{2413}{2}$ is
represented by $2-41-3$. We stick to our notation so that it can be generalized
for $d$-permutations.

\begin{figure}[htb]
	\center{\resizebox{0.5\textwidth}{!}{\input{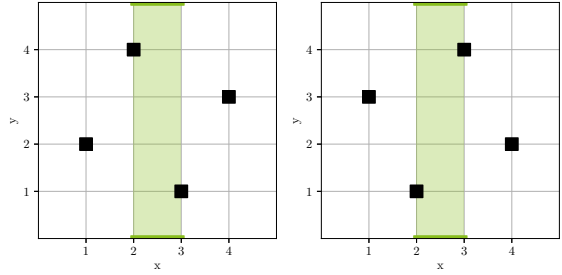}}}
	\caption{Baxter permutation forbidden vincular patterns: $\vinpat{2413}{2}$
	and $\vinpat{3142}{2}$. The adjacency is materialized by a vertical (green)
	strip.}\label{fig:baxter-patterns}
\end{figure}

Baxter permutations (introduced by Glen Baxter~\cite{Bax:65}) are exactly permutation that avoids $\vinpat{2413}{2}$ and $\vinpat{3142}{2}$ (see Figure~\ref{fig:baxter-patterns}):

$$B_n:=S_n(\vinpat{2413}{2},\vinpat{3142}{2}).$$

$$|B_n| = \sum_{k=1}^{n}\frac{\binom{n+1}{k-1}\binom{n+1}{k}\binom{n+1}{k+1}}{\binom{n+1}{1}\binom{n+1}{2}}.$$

The first terms of $(B_n)$ are 1, 2, 6, 22, 92, 422, 2074 (sequence A001181 in
OEIS).

\begin{figure}[htb]
	\center{\resizebox{0.5\textwidth}{!}{\input{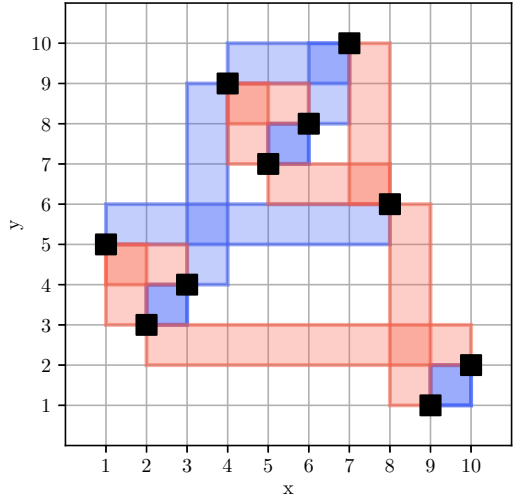}}}
	\caption{Example of a Baxter permutation. At each ascent (resp. descent) we
		associate a blue (resp. red) vertical rectangle, called \emph{slice}, and
		we associate a blue (resp. red) horizontal rectangle to each ascent (resp.
		descent) of the inverse permutation.}\label{fig:baxter2d}
\end{figure}

Figures~\ref{fig:baxter2d} and the 2 first permutations of Figure~\ref{fig:2dsep} give examples of Baxter permutations.

\subsection{Separable $d$-permutations}

A \emph{$d$-direction} (or simply a \emph{direction}) $\bdir$ is a word on the
alphabet $\{+,-\}$ of length $d$. A direction is \emph{positive} if its first
entry is positive.

Let $\bsig$ and $\bpi$ two $d$-permutations and $\bdir$ a positive direction.
The \emph{$d$-sum} with respect to direction $\bdir$ is the $d$-permutation:
$$\bsig \oplus^\bdir \bpi := \basig_2 \oplus^\bdir_2 \bapi_2,\dots, \basig_{d}
	\oplus^\bdir_d \bapi_{d},$$
where $\oplus^\bdir_i$ is $\oplus$ if $\bdir_i=+$ and $\ominus$ if $\bdir_i=-$.

\begin{figure}[htb]
	\centering{\resizebox{0.7\textwidth}{!}{\input{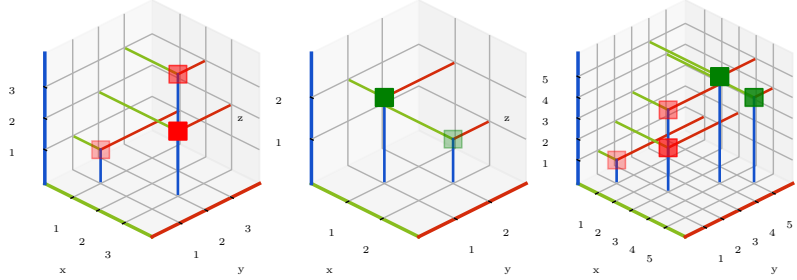}}}
	\caption{On the left a permutation $p_1 = \perm{132}{132}$ and in the middle a permutation $p_2 = \perm{12}{21}$.
		$p_1$ and $p_2$ are separable $3$-permutations because $p_1 = \perm{1}{1} \oplus^{(+++)} (\perm{1}{1} \oplus^{(+--)} \perm{1}{1})$ and $p_2 = \perm{1}{1} \oplus^{(+-+)} \perm{1}{1}$.
		On the right, their $d$-sum with respect to direction $(+++)$ is $\perm{132}{132} \oplus^{(+++)} \perm{21}{21} = \perm{13254}{13254}$ which is still separable.}
	\label{fig:3d-sep-sum}
\end{figure}

A \emph{separable $d$-permutation} is a $d$ permutation of size 1 or the $d$-sum of two separable $d$-permutations.
These definitions are illustrated in Figure~\ref{fig:3d-sep-sum}.

\newpage

As we have seen previously, for $d=2$ every permutation of size at most $3$ is
separable and these permutations are characterized by the avoidance of 2
patterns of size 4. For $d=3$, it’s no longer true that all 3-permutations of
size 3 are separable. The 8 3-permutations of size 3 that are not separable are
$\patt$ (see Figure~\ref{fig:patt3}). In fact, these 8 permutations together
with the 2 patterns of length 4 characterize exactly separable $d$-permutations
for any $d\geq 3$ as shown in~\cite{asinowski2010separable}. We restart their
result with our formalism:
\begin{theorem}~\cite{asinowski2010separable}
	Let $Sep_n^{d-1}$ be set of separable $d$-permutations of size $n$. $$Sep_n^{d-1}=S_n^{d-1}(\patt, 2413, 3142).$$
\end{theorem}

Explicit enumeration formulas have been proposed in~\cite{asinowski2010separable}:
$$|Sep_n^{d-1}| = \frac{1}{n-1} \sum_{k=0}^{n-2}\binom{n-1}{k}\binom{n-1}{k+1}(2^{d-1}-1)^k(2^{d-1})^{n-k-1}.$$

\begin{table}[!htb]
	\center{\begin{tabular}{|c||c|c|c|c|}
			\hline
			n / d & 2    & 3      & 4        & 5          \\
			\hline
			\hline
			1     & 1    & 1      & 1        & 1          \\
			\hline
			2     & 2    & 4      & 8        & 16         \\
			\hline
			3     & 6    & 28     & 120      & 496        \\
			\hline
			4     & 22   & 244    & 2248     & 19216      \\
			\hline
			5     & 90   & 2380   & 47160    & 833776     \\
			\hline
			6     & 394  & 24868  & 1059976  & 38760976   \\
			\hline
			7     & 1806 & 272188 & 24958200 & 1887736816 \\
			\hline
		\end{tabular}}
	\caption{Values of $|Sep_n^{d-1}|$ for the first values of $n$ and $d$.}
\end{table}

Now we give a new characterization of separable $d$-permutations (Theorem~\ref{thm:cara-sep}).
Thanks to this result, checking if a $d$-permutation is separable is simpler:
we only need to check if it avoids the dimension 3 patterns and then if it
avoids the dimension 2 patterns only on $d-1$ projections instead of
$(d-1)\times(d-2)/2$ projections.

\begin{theorem}\label{thm:cara-sep}
	$$Sep_n^{d-1} = S_n(2413, 3142)^{d-1} \cap S_n^{d-1}(\patt).$$
\end{theorem}

\begin{proof}
	To show this result, we only need to show that for any $\bsig \in S_n^{d-1}(\patt)$ and any $1<i<j\leq n $, if $\proj_{i,j}(\bsig)$ contains one of the patterns $2413, 3142$ then $\sigma_j$ also does.

	So let $\bsig\in S_n^{d-1}(\patt)$ and $1 <i,j\leq n$ such that
	$\proj_{i,j}(\bsig)$ contains the pattern $2413$ (the other case being
	identical). Let $p_1, p_2, p_3, p_4\in P_\bsig$  an occurence of this
	pattern such that  $x(p_1) < x(p_2) < x(p_3) < x(p_4)$. The projection of
	$p_1$ and $p_2$ in the plane $(x_i,x_j)$ are in the same quadrant of the
	projection of $p_3$. Since $\bsig$ avoids $\patt$ and by
	Remark~\ref{rem:patt}, $x(p_3)$ is not between $x(p_1)$ and $x(p_2)$.

	Applying the same argument on the 3 other triplets of points, we get that $x(p_1)$ is not between $x(p_2)$ and $x(p_4)$, $x(p_3)$ is not between $x(p_1)$ and $x(p_2)$, $x(p_4)$ is not between $x(p_1)$ and $x(p_3)$.

	There is only two orders that satisfy these four constrains: $x(p_1)< x(p_2) < x(p_3) < x(p_4)$ and $x(p_4) < x(p_3) < x(p_2) < x(p_1)$. In the first case, the four points induce the pattern $2413$ on $\proj_{1,j}$ in the second case it is the pattern $3142$.

	Hence if $\proj_{i,j}(\bsig)$ contains a forbidden pattern so does $\proj_{1,j}(\bsig)=\sigma_j$.
\end{proof}

\subsection{Baxter $d$-permutations}
We now generalize Baxter permutations to larger dimensions. To do so we
introduce some formalism that will ease the definition of these
$d$-permutations.

Given $P_\bsig$ the diagram of a $d$-permutation $\bsig$, two points $p_i,p_j$
of
$P_\bsig$ are \emph{$k$-adjacent} if they differ by one on their $k$-th
coordinate, and $k$ is said to be the \emph{type} of the adjacency.
The \emph{direction} of $p_i,p_j$ is the signs of the coordinates of $p_j-p_i$.
Given two adjacent points $p_i$ and $p_j$, the \emph{slice} of $p_i,p_j$ is the
$d$-dimensional box  with $p_i$ and $p_j$ as corners. a slice $p_i,p_j$ is of
\emph{type} $k$ is $p_i,p_j$ are $k$-adjacent. The \emph{direction} of a slice
$p_i,p_j$ is the direction of $p_i,p_j$ if $x(p_i)<x(p_j)$ and the direction of
$p_j,p_i$ otherwise. A \emph{cell} is a unitary cube whose corners have
coordinates are integers. A single slice can have multiple types.
For instance, if a slice is a cell it is of all possible types.

For $d=2$, an ascent in a permutation corresponds to an adjacency of type 1
(which corresponds to the $x$-axis) with direction $++$, a descent is an
adjacency of type 1 with direction $+-$. An adjacency of type 2 (which
corresponds to the $y$-axis) with direction $+-$ corresponds to an ascent in
the inverse permutation.

In Figure~\ref{fig:baxter2d}, slides of direction $++$ are represented in blue and those of type $+-$ in blue.

\begin{definition}
	A $d$-permutation is \emph{well-sliced} if each slice intersects exactly one slice of each type and two intersecting slices share the same direction.
\end{definition}

One can observe that the Baxter permutation of Figure~\ref{fig:baxter2d} is
well-sliced.

\begin{definition}
	A \emph{Baxter $d$-permutation} is a $d$-permutation such that each of its $d' \leq d$ projection is well-sliced.
\end{definition}

\begin{figure}[htb]
	\center{\resizebox{0.8\textwidth}{!}{\input{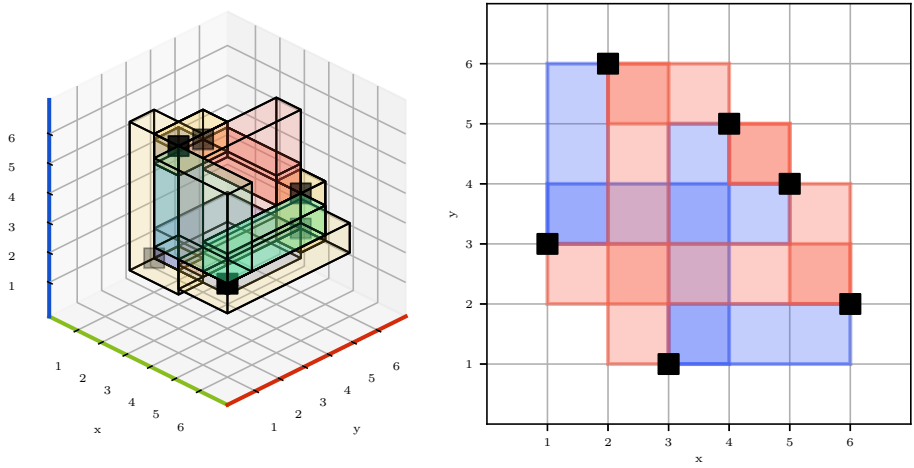}}}
	\caption{On the left, \perm{342651}{156243}, an example of 3-permutation that is
		well-sliced but not Baxter since its projection on the plane (y,z) ({\color{axeX} 361542}) on the right is not well-sliced.}\label{fig:notVeryWellSliced}
\end{figure}

By definition, if $d$-permutation is Baxter it is also the case for all its projections of smaller dimensions. On the other hand, a $d$-permutation can be well-sliced and have projections that are not well-sliced. Take for instance the $3$-permutation \perm{342651}{156243}. Its projection on the plane (y,z) is  {\color{axeX} 361542} is not well-sliced since it is not a Baxter permutation. (see Figure~\ref{fig:notVeryWellSliced})

Table~\ref{tab:seqBax} gives the first values of $|B_n^{d-1}|$.
\begin{table}[!htb]
	\center{\begin{tabular}{|c||c|c|c|c|}
			\hline
			n / d & 2    & 3      & 4     & 5       \\
			\hline
			\hline
			1     & 1    & 1      & 1     & 1       \\
			\hline
			2     & 2    & 4      & 8     & 16      \\
			\hline
			3     & 6    & 28     & 120   & 496     \\
			\hline
			4     & 22   & 260    & 2440  & 20816   \\
			\hline
			5     & 92   & 2872   & 59312 & 1035616 \\
			\hline
			6     & 422  & 35620  &       &         \\
			\hline
			7     & 2074 & 479508 &       &         \\
			\hline
		\end{tabular}}
	\caption{Values of $|B_n^{d-1}|$ for the first values of $n$ and
	$d$.}\label{tab:seqBax}
\end{table}

\begin{figure}[htb]
	\center{\input{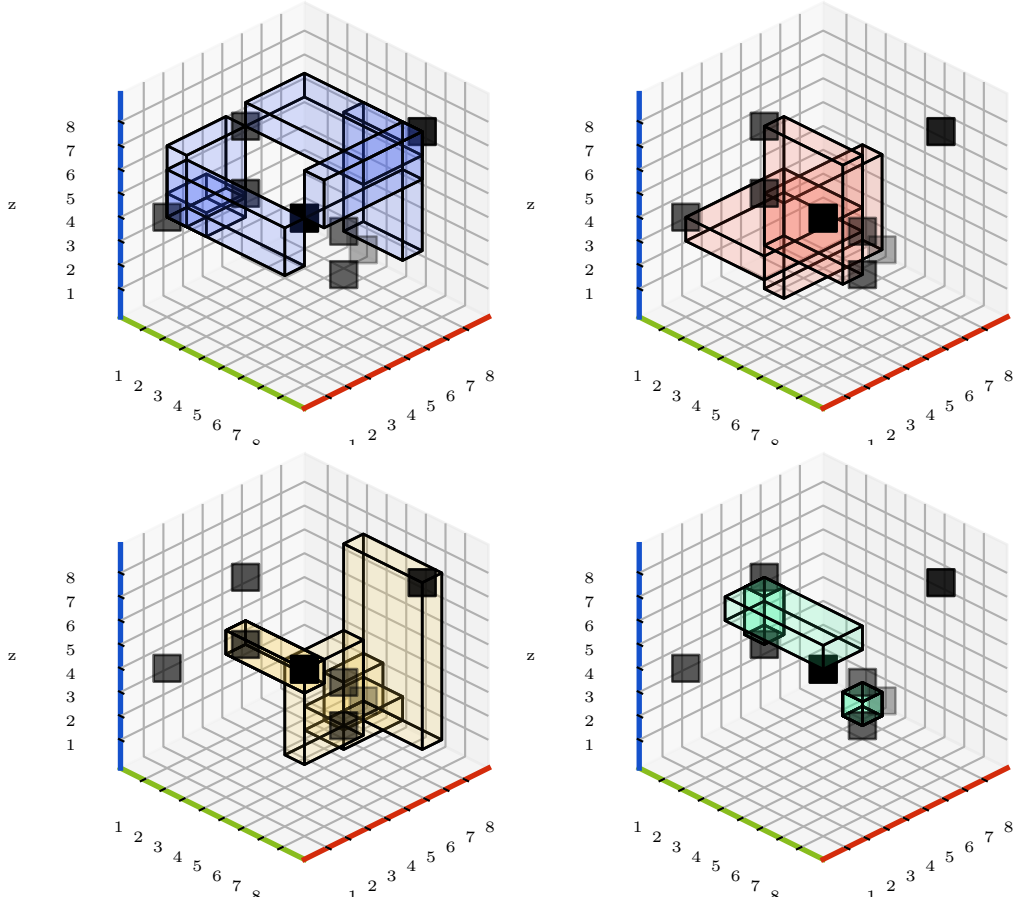}}
	\caption{\permbaxter: an example of a Baxter 3-permutation, together with
		its slices of different types.}\label{fig:paral1}
\end{figure}


In order to characterize Baxter $d$-permutation, let us introduce generalized
vincular patterns.

\begin{definition}
	A \emph {generalized vincular pattern} \vinpat{\bpi}{X_1,\cdots,X_d} is a
	permutation $\bpi$ together with a list of subsets of $[k-1]$ $X_1,\cdots,X_d$
	called \emph{adjacencies}.
	Given $\bsig$ a $d$-permutation, we say that $p_1,\cdots p_k \in P_\bsig$ is an
	\emph{occurrence} of the pattern
	\vinpat{\bpi}{X_1,\cdots,X_d} if $p_1,\cdots p_k$ is an occurrence of $\pi$
	and  if it satisfies the adjacency constrains: for each $k$ and each $i\in
		X_k$: the $i$-th and $(i+1)$-th points with respect to the order along the axis
	$k$ are $k$-adjacent. We say that $\bsig$ a $d$-permutation contains the pattern
	\vinpat{\bpi}{X_1,\cdots,X_d'} (of dimension $d'$), if at least one direct
	projection of
	dimension $d'$ of $\bsig$ contains an occurrence of the pattern
	\vinpat{\bpi}{X_1,\cdots,X_d'}.
\end{definition}

\begin{figure}[htb]
	\center{\resizebox{0.5\textwidth}{!}{\input{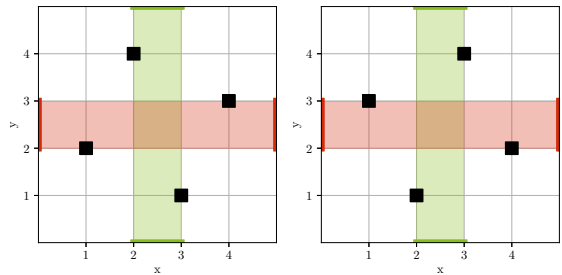}}}
	\caption{Baxter permutation can also be characterized by these two
		generalized vincular forbidden patterns:
		$\vinpatb{2413}{2}{2}$ and
		$\vinpatb{3142}{2}{2}$.}\label{fig:baxter-patterns-bis}
\end{figure}

It is well known that
$S_n(\vinpat{2413}{2}) = S_n(\vinpatb{2413}{2}{2})$  and $S_n(\vinpat{3142}{2})
	= S_n(\vinpatb{3142}{2}{2})$ (see
Figure~\ref{fig:baxter-patterns-bis}). Every occurrence of
$\vinpatb{2413}{2}{2}$
is clearly an occurrence of \vinpat{2413}{2}. The converse
is obtained thanks to the following observation: if $i_1,i_2,i_3,i_4$ is an
occurrence of $\vinpat{2413}{2}$ in $\sigma$, let $i'_1$ such that $i'_1<i_2$
and $\sigma(i_1) \leq \sigma(i'_1) < \sigma(i_4)$, such that  $\sigma(i'_1)$ is
maximal.
Let $i'_4 = \sigma^{-1}(\sigma(i'_1)+1)$. We have that $i'_1,i_2,i_3,i'_4$ is
an
occurrence of $\vinpatb{2413}{2}{2}$.

It follows that
$$B_n=S_n(\vinpatb{2413}{2}{2},\vinpatb{3142}{2}{2}).$$

As a warm-up for the rest of the
section, let us reprove that our definition of Baxter $d$-permutations coincides with the classical one.

\begin{proposition}\label{pro:baxter-well-sliced}
	A permutation is a Baxter permutation if and only if it is well-sliced.
\end{proposition}
\begin{proof}
	As shown above $B_n=S_n(\vinpatb{2413}{2}{2},\vinpatb{3142}{2}{2}).$
	If a permutation contains one of the above patterns, then it contains 2
	intersecting slices of different directions, hence it is not  well-sliced.
	Now let consider a permutation $\sigma$ that is not well-sliced and let us show
	that it contains a forbidden pattern.
	As it is not well-sliced it contains {\bf (i)} a pair of intersecting slices of
	different directions, {\bf (ii)} it contains a slice that intersects two other
	slices
	or {\bf (iii)} it contains a slice that doesn't intersect any other slices.

		{\bf (i)}: Any occurrence of two slices of different directions is an
	occurrence of
	one of the two forbidden patterns.

		{\bf (ii)}: Let $p_1,p'_1, p_2, p_3,p_4,p'_4$ such that $p_2,p_3$ is a vertical
	slice, and
	$p_1,p_4$ and $p'_1$, $p'_4$ are two horizontal slices intersecting the slice
	$p_2,p_3$.
	Since we have treated the case {\it (i)} we can assume that the 3
	slices are of the same type and without loss of generality we can assume that
	this type is $(++)$. Observe that $p_1,p'_1, p_4, p'_4$ are four different
	points but this set of points may intersect the point set $\{p_2,p_3\}$.
	Nevertheless we can assume that $p_1$ and $p'_1$ are on the left of $p_3$ and
	$p_4$ and $p'_4$ are on the right of $p_2$. We can
	also assume without loss of generality that $p'_1$ and $p'_4$ are below $p_1$
	and $p_4$. Hence $p_1,p_2,p_3,p'_4$ are 4 different points and we can then
	observe that this point set is an occurrence of
	$\vinpat{3142}{2}$, hence $\sigma$ contains $\vinpatb{3142}{2}{2}$.

		{\bf (iii)}: Let us show this case can't occur. In other words, let us show
	that every vertical slice intersects at least one horizontal slice. Without
	loss of generality let us just consider the case of an
	ascent. Let $i_1$ such that $\sigma(i_1)<\sigma(i_1+1)$. Let $i_2$ such that
	$i_2 \leq i_1$ such that $\sigma(i_1)\leq \sigma(i_2) < \sigma(i_1+1)$ and such
	that $\sigma(i_2)$ is maximal. Let $i_3=\sigma^{-1}(\sigma(i_2)+1)$. By
	construction $i_3 \geq i_2$. Hence, the vertical slice $p_{i_1},p_{i_1+1}$
	intersects the horizontal slice $p_{i_2},p_{i_3}$. Contradiction.
\end{proof}

The action of the symmetries of the hypercube extends naturally on the
generalized vincular patterns. We can remark that $\sbaxpa =
	\{\vinpatb{2413}{2}{2},\vinpatb{3142}{2}{2}\}$, hence, $B_n=S_n(\sbaxpa)$.

\begin{figure}[htb]
	\center{\resizebox{1\textwidth}{!}{\input{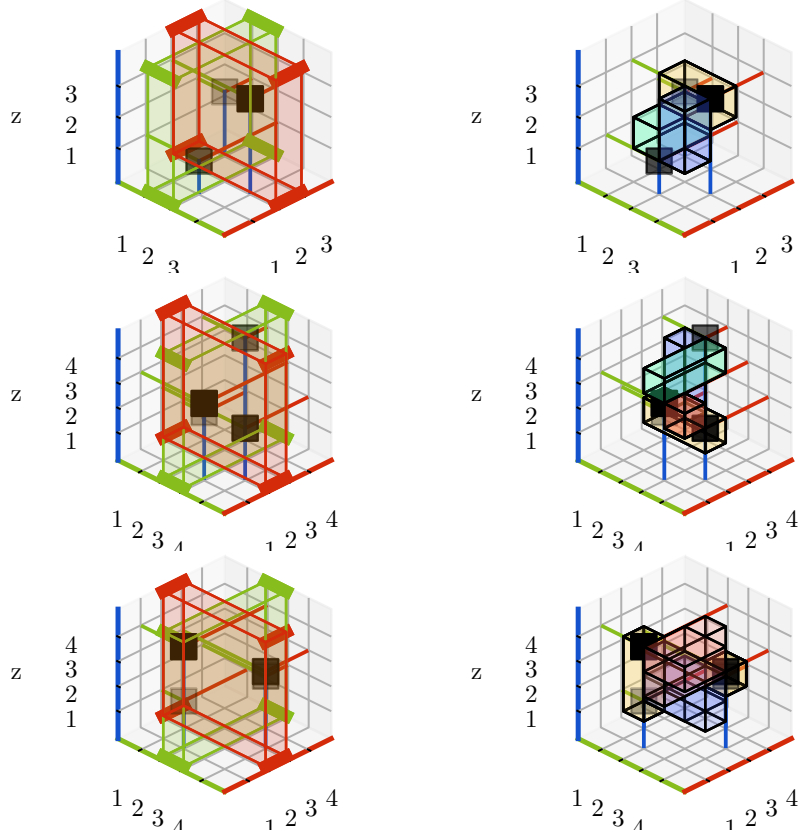}}}
	\caption{On the left, the 3 3-dimensional vincular pattern forbidden in
		Baxter
		$d$-permutations: \baxpb,\baxpc,\baxpd. The adjacency constrains are
		materialized by
		boxes orthogonal to the concerned axes.  On the right the
		corresponding 3-permutations
		with all its slices. One can observe that it is not
		well-sliced because the two first have a lack of slice
		intersections and the last one a bad
		intersection.}\label{fig:baxter3d-patterns}
\end{figure}

\begin{theorem}\label{thm:cara-bax}
	$$B_n^{d-1} = S_n^{d-1}(\sbaxpa, \sbaxpb, \sbaxpc, \sbaxpd).$$
\end{theorem}
Figure~\ref{fig:baxter3d-patterns} depicts on occurrence of each class of forbidden patterns of dimension 3.
The list of all symmetries of these patterns are given in the Appendix~\ref{sec:sym_bx}.

\begin{proof}
	Let us start with the easy inclusion:

	$\subseteq$: Let $\bsig$ a $d$-permutation that contains one of the forbidden patterns.
	If a $d$-permutation contains one of the forbidden patterns \sbaxpa (resp. \sbaxpd), then at least one of its 2 dimensional (resp. 3 dimensional) projection is not well sliced since these patterns are witnesses of the intersections of 2 slices of different directions.
	Hence $\bsig$ isn't Baxter.

	If $p_1,p_2,p_3$ (resp. $p_1,p_2,p_3,p_4$) is an occurrence of the pattern \baxpb  (resp.\baxpc) on one of the 3-dimensional projection of $\bsig :=\bsig_{3}$, then
	the slices $p_1,p_2$ and $p_1,p_3$ (resp. $p_1,p_4$ and $p_2,p_3$) doesn't intersect. We remark that on $\proj_{x,y}(\bsig_3)$ the corresponding slices intersect. Hence, either there is no other intersection of the slices $p_1,p_2$ (resp. $p_1,p_4$) in $\bsig_{3}$ and $\bsig_{3}$ is not well sliced, either slice intersects another slice in $\bsig_{3}$ and in this case  the slice $p_1,p_2$ (resp. $p_1,p_4$) intersects two slices in $\proj_{x,y}(\bsig_3)$. In both cases $\bsig$ is not Baxter.
	We can apply the same reasoning of all symmetries of \baxpb and \baxpc.
	Now let us consider the other inclusion.

	$\supseteq$: Let $\bsig$ a $d$-permutation that is not Baxter and let us prove that it contains one of the forbidden patterns. Let us consider the three following sub-cases:

	\begin{itemize}
		\item {\bf (i) there are 2 intersecting slices of different directions
		.}
		      Assume without loss of generality that the slice $p_2,p_3$ of
		      type $x$ intersects slice the $p_1,p_4$ of type $y$. If the signs
		      of
		      the direction of the slices are different on $x$ or $y$ then
		      $p_1,
		      p_2, p_3, p_4$ is an occurrence a forbidden pattern in $\sbaxpa$
		      in $\proj_{xy}(\bsig)$.
		      So now let us assume that the directions of these two slices
		      share the same signs on coordinates $x$ and $y$ but differ on a
		      third coordinate. Without loss of generality assume that the
		      third coordinate is $z$ and in $\proj_{xyz}(\bsig)$ the direction
		      for the first one is $(+++)$ and $(++-)$ for the second.
		      First observe that since these two slices intersect each other and are of different types, $p_1,p_2,p_3,p_4$ are 4 different points and we have $x(p_1) < x(p_2) < x(p_3) < x(p_4)$ and $y(p_2) < y(p_1) < y(p_4) <y(p_3)$. Moreover we have that $z(p_2) <z(p_3)$ and $z(p_4) < z(p_1)$. If $z(p_1)$ and $z(p_4)$ are between $z(p_2)$ and $z(p_3)$ then $\proj_{xz}(\bsig)$ contains a forbidden pattern in $\sbaxpa$.
		      If $z(p_2)$ and $z(p_3)$ are between $z(p_4)$ and $z(p_1)$ then $\proj_{yz}(\bsig)$ contains a forbidden pattern in $\sbaxpa$. If it is not the case, then we have $z(p_2) < z(p_4) < z(p_3) < z(p_1)$ or $z(p_4) < z(p_2) < z(p_1) < z(p_3)$. In these last two cases $p_1,p_2,p_3,p_4$ is an occurrence of a forbidden pattern of $\sbaxpd$ in $\proj_{xyz}(\bsig)$.

		\item {\bf (ii) there is a slice that intersects 2 slices of the same
		type.}
		      Assume that there is a slice $p_1,p_6$ of type $y$ that intersect
		      2 slices of type $x$, $p_2,p_3$ and $p_4,p_5$, such that $x(p_1)<
		      x(p_2) \cdots < x(p_6)$. Since we have already treated the case
		      of intersections of different directions, we can assume that
		      these 3 slices share the same direction and without loss of
		      generality we can assume that this is the positive direction.
		      This implies that $y(p_3),y(p_5) > y(p_6)$ and $y(p_2),y(p_4) <
		      y(p_1)$. Hence $p_1,p_3,p_4,p_6$ is an occurrence of
		      \vinpatb{3142}{.}{2} in $\proj_{xy}(\bsig)$. Hence $\bsig$
		      contains a pattern of $\sbaxpa$.

		\item {\bf (iii) there is a slice that intersects no slice of a given
		type.}
		      Without loss of generality let us consider the all positive direction. Assume there is a $x$-slice $(p_2,p_3)$ that doesn't intersect any $y$-slice.  Let us consider $\proj_{xy}(\bsig)$. If $\bsig$ is not Baxter, $\proj_{xy}(\bsig)$ contains a forbidden pattern $\sbaxpa$. Otherwise, in $\proj_{xy}(\bsig)$, the slice $(p_2,p_3)$ intersects exactly one slice. Let $p_2,p_3$ such that the slice $(p_1,p_4)$ intersects the slice $(p_2,p_3)$ in $\proj_{xy}(\bsig)$. Remark that the point $p_1$ may be equal to $p_2$. Since these two slices doesn't intersect in $\bsig$, there must be a third coordinate, $z$ without loss of generality, such that $z(p_1),z(p_4) \leq z(p_2)$ or $z(p_1),z(p_4) > z(p_3)$.
		      If $p_1=p_3$, then the 3 points form an occurrence of a forbidden pattern in $\sbaxpb$ otherwise the 4 points form an occurrence of a forbidden pattern in $\sbaxpc$.
	\end{itemize}
\end{proof}

As all the patterns involved in the previous theorem are of dimension 2 or 3, we get the following corollary:
\begin{corollary}
	A $d$-permutation is Baxter if and only if all its projections of dimensions 2 or 3 are well-sliced.
\end{corollary}

\subsection{Anti and complete Baxter $d$-permutations}
\begin{figure}
	\center{\resizebox{0.4\textwidth}{!}{\input{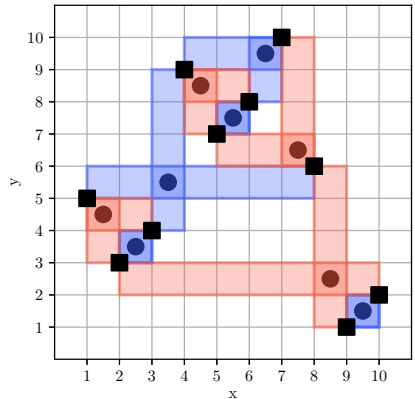}}}
	\caption{The Baxter permutation $5\, 3\, 4\, 9\, 7\, 8\, 10\, 6\, 1\, 2$
		(square points) together with its associate anti-Baxter permutation
		(circle
		points ) $4\, 3\, 5\, 8\, 7\, 9\, 6\, 2\, 1$. The corresponding the
		complete
		Baxter permutation (all points together) is
		$9\, 8\, 5\, 6\, 7\, 10\, 17\, 16\, 13\, 14\, 15\, 18\, 19\, 12\, 11\,
			4\, 1\,
			2\, 3$.}\label{fig:anti2dV}
\end{figure}

In a Baxter permutation $\sigma$, each vertical slice intersects exactly one horizontal
slice, these intersections are cells (square of width 1) (See for instance
Figure~\ref{fig:anti2dV}). Let $P'_\sigma$ the set of centers of these cells.
If
we combine $P_\sigma$ and $P'_\sigma$ we obtain the diagram of permutation of
size
$2n+1$ (on a finer grid). These permutations are often called \emph{complete}
Baxter permutations and was introduced by Baxter and
Joichi~\cite{baxterjoichi} under the name \emph{$w$-admissible}
permutations. What we call here Baxter permutations are sometimes called
\emph{reduced} Baxter permutations.

\begin{figure}[htb]
	\center{\resizebox{0.5\textwidth}{!}{\input{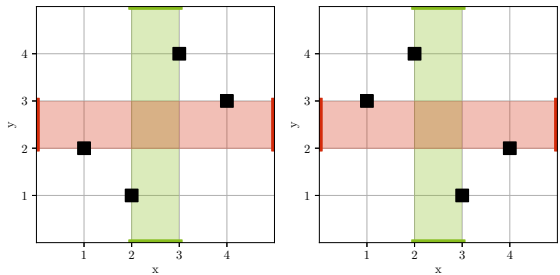}}}
	\caption{Forbidden patterns in anti-Baxter permutations:
		$\vinpatb{2143}{2}{2}$ and
		$\vinpatb{3412}{2}{2}$.}\label{fig:anti-baxter-patterns-bis}
\end{figure}

The permutations corresponding to $P'_\sigma$ are called \emph{anti-Baxter}
permutations. These permutations are exactly the ones avoiding
$\vinpatb{2143}{2}{.}$ and $\vinpatb{3412}{2}{.}$ as shown
in~\cite{asinowski2013orders}. As for Baxter patterns:
$S_n(\vinpatb{2143}{2}{.},\vinpatb{3412}{2}{.}) =
	S_n(\vinpatb{2143}{2}{2},\vinpatb{3412}{2}{2})$ (see \cite[Lemma
	3.5]{asinowski2013orders} and Figure~\ref{fig:anti-baxter-patterns-bis}). The
enumeration of this class of permutation has
been given in~\cite{asinowski2013orders}

Let us generalize these definitions of anti-Baxter and complete Baxter to
higher dimensions. For that purpose let us start with the following property.

\begin{proposition}\label{pro:intersect}
	Let $\bsig$ a well-sliced $d$-permutation. Given a slice $p_1,p'_1$ of type
	$1$, let
	$(p_i,p'_i)$ be the slices of type $i \in [d]$ that intersect $p_1,p'_1$. The
	intersection of all these slices is the cell $q,q'$, where
	$x_i(q):=x_i(p_i)$ and $x_i(q'):=x_i(p'_i)$.
\end{proposition}
\begin{proof}
	First observe that the cell $q,q'$ is included in each slice $p_i,p'_i$. Hence  the cell $q, q'$ is included in the intersection of all slices $p_i, p'_i$.

	Since every slice $p_j,p'_j$ intersects the slice $p_i,p'i$, we have
	$\max(\min(x_i(p_i),x_i(p'_i)),\min(x_i(p_j),x_i(p'_j))) <
		\min(\max(x_i(p_i),x_i(p'_i)),\max(x_i(p_j),x_i(p'_j)))$. Moreover, since the
	slice $p_i,p'_i$ is of width 1 with respect to axis $i$ and all the other have
	a width greater or equal to one, we have
	$\min(x_i(p_j),x_i(p'_j)) \leq \min(x_i(p_i),x_i(p'_i))$ and
	$\max(x_i(p_j),x_i(p'_j)) \geq \max(x_i(p_i),x_i(p'_i))$.
	Hence the intersection of the projections of slices on the axis $i$ is the
	interval $[\min(x_i(p_i),x_i(p'_i)), \max(x_i(p_i),x_i(p'_i))]$. Hence the
	intersection of the considered slices is included in the slice $q,q'$.
\end{proof}

\begin{figure}
	\center{\resizebox{1\textwidth}{!}{\input{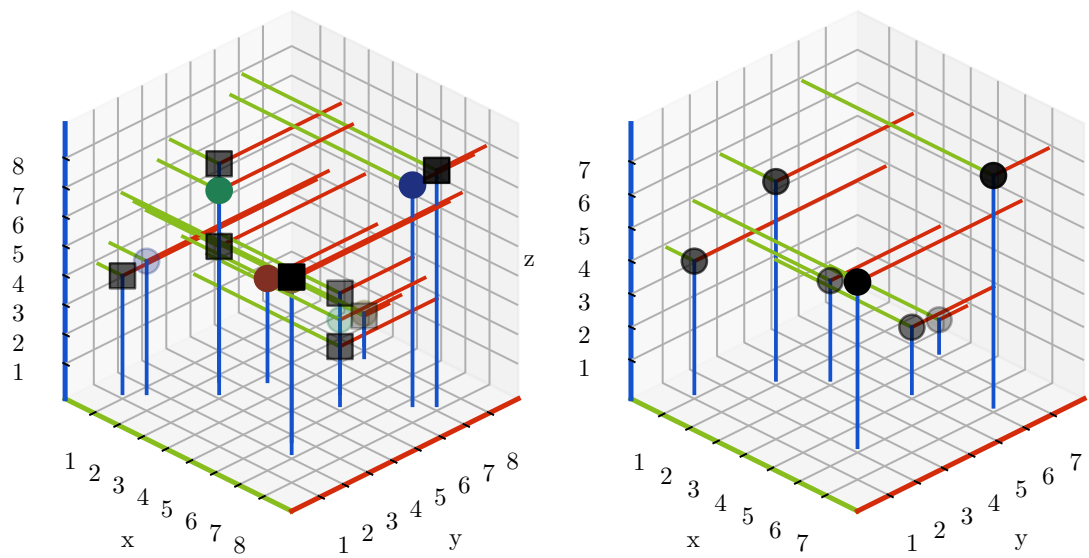}}}
	\caption{On the left, the complete Baxter 3-permutation
		$\perm{14386527}{47513268}$
		with its cell
		(circle) points. Each cell point corresponds to the triple intersection
		of slices of the same type (see
		Figure~\ref{fig:paral1}). On the right, the
		anti-Baxter 3-permutation $\perm{1347526}{4631257}$ associate
		with the Baxter permutation of Figure~\ref{fig:bax3dV}.}\label{fig:bax3dV}
\end{figure}

Given a Baxter $d$-permutation $\bsig$, for every slice of type $1$, let us
associated the \emph{intersecting cell} defined from
Property~\ref{pro:intersect} (see Figure~\ref{fig:bax3dV}).
Let $P'_\bsig$ be the set of centers of intersecting cells. Since every slice
of any type contains exactly one intersecting cell, $P'_\bsig$ defines a
$d$-permutation, and we call the $d$-permutations obtained this way
\emph{anti-Baxter} $d$-permutations (see Figure~\ref{fig:bax3dV}). Again, this
definition coincides with the classical one. If we combine $P_\sigma$ and
$P'_\sigma$ we obtain the diagram of $d$-permutation of size $2n+1$ (on a finer
grid). We naturally call these $d$-permutations  \emph{complete Baxter
$d$-permutations}.

As for Baxter
$d$-permutations, a projection of an anti-Baxter (resp. a complete
Baxter)
$d$-permutation is also an anti-Baxter (resp. a complete Baxter)
$d'$-permutation. We denote by
$A_n^{d-1}$ the set of anti-Baxter $d$-permutations of size $n$. The first
values of $A_n^{d-1}$ are given in Table~\ref{tab:seqAnti}.

\begin{table}[!htb]
	\center{\begin{tabular}{|c||c|c|c|c|}
			\hline
			n / d & 2    & 3     & 4    & 5    \\
			\hline
			\hline
			1     & 1    & 1     & 1    & 1    \\
			\hline
			2     & 2    & 4     & 8    & 16   \\
			\hline
			3     & 6    & 36    & 216  & 1296 \\
			\hline
			4     & 22   & 444   & 7096 &      \\
			\hline
			5     & 88   & 5344  &      &      \\
			\hline
			6     & 374  & 64460 &      &      \\
			\hline
			7     & 1668 &       &      &      \\
			\hline
		\end{tabular}}
	\caption{Values of $|A_n^{d-1}|$ for the first values of $n$ and
	$d$.}\label{tab:seqAnti}
\end{table}


%


\section{Conclusion and perspectives}\label{sec:conclusion}
In this paper we have started to consider pattern avoidance in $d$-permutations
and we have generalized Baxter permutations in that context. These first steps
give rise to a large number of open problems, some probably hard but some probably
very tractable.

The enumeration of $d$-permutations avoiding the smallest patterns is widely
open, starting from the smallest one: $\perm{12}{12}$. Moreover, as presented
many known enumeration sequences seems to match with several permutations
family. Clearly, there are several bijections to find to explain these
sequences.

Around Baxter $d$-permutations, it is a large field of research that is open.
Let us mention few examples of questions related to this family of permutations.  Clearly, the first-expected result would be the enumeration of Baxter $d$-permutations. As mention in the introduction, Baxter permutations are in bijection we several interesting combinatorial objects. A very natural question would be:  which of these bijections can be extended to $d$-Baxter permutations.
For instance, Baxter permutations are in bijection with boxed arrangements of
axis-parallel segments in $\mathbb{R}^2$~\cite{felsner2011bijections}.
In~\cite{felsner2020plattenbauten}, they studied boxed arrangements of
axis-parallel segments in $\mathbb{R}^3$. Are there some links between Baxter
$d$-permutations boxed arrangements in $\mathbb{R}^{2^{d-1}}$?

We were able to characterize Baxter $d$-permutations with forbidden vincular patterns. This question remains open for anti-Baxter $d$-permutations.

In addition to that, several classes related to Baxter permutations have
received some attention: \emph{doubly alternating} Baxter
permutations~\cite{guibert2000doubly}, Baxter
\emph{involutions}~\cite{fusy2012bijective}, \emph{semi} and \emph{strong} Baxter
permutations~\cite{bouvel2018semi}, \emph{twisted} Baxter permutations~\cite{west2006enumeration}. Once again, can some of these classes
be extended and enumerated in higher dimension?

We have developed a module base on Sage to work with $d$-permutation~\url{https://plmlab.math.cnrs.fr/bonichon/multipermutation}. We hope this tool will help the community to investigate the above problems.

\section*{Acknowledgment} The authors would like to thank Eric Fusy, Valentin Feray, Mathilde Bouvel and Olivier Guibert for fruitful discussions. This work was partially supported by ANR grant 3DMaps ANR-20-CE48-0018.

\bibliography{bibliography}{}

\begin{thebibliography}{10}

\bibitem{aldred2005permuting}
Robert~EL Aldred, Mike~D Atkinson, Hans~P van Ditmarsch, Chris~C Handley,
  Derek~A Holton, and DJ~McCaughan, Permuting machines and priority queues,
  {\em Theoretical computer science} {\bf 349}(3) (2005), 309--317.

\bibitem{asinowski2013orders}
Andrei Asinowski, Gill Barequet, Mireille Bousquet-M{\'e}lou, Toufik Mansour,
  and Ron~Y Pinter, Orders induced by segments in floorplans and $(2-14-3,
  3-41-2) $-avoiding permutations, {\em The Electronic Journal of
  Combinatorics} {\bf 20}(2) (2013), P35.

\bibitem{asinowski2010separable}
Andrei Asinowski and Toufik Mansour, Separable d-permutations and guillotine
  partitions, {\em Annals of Combinatorics} {\bf 14}(1) (2010), 17--43.

\bibitem{atkinson1995priority}
MD~Atkinson, SA~Linton, and LA~Walker, Priority queues and multisets, {\em the
  electronic journal of combinatorics}  (1995), R24--R24.

\bibitem{atkinson1993permutational}
MD~Atkinson and Murali Thiyagarajah, The permutational power of a priority
  queue, {\em BIT Numerical Mathematics} {\bf 33}(1) (1993), 1--6.

\bibitem{aval2021baxter}
Jean-Christophe Aval, Adrien Boussicault, Mathilde Bouvel, Olivier Guibert, and
  Matteo Silimbani, Baxter tree-like tableaux, 2021.

\bibitem{Bax:65}
G.~Baxter, On fixed points of the composite of commuting functions, {\em
  Proceedings of the American Mathematical Society} {\bf 15} (1964), 851--855.

\bibitem{baxterjoichi}
Glen Baxter and JT~Joichi, On permutations induced by commuting functions, and
  an embedding question, {\em Mathematica Scandinavica} {\bf 13}(2) (1963),
  140--150.

\bibitem{bonichon2010baxter}
Nicolas Bonichon, Mireille Bousquet-M{\'e}lou, and {\'E}ric Fusy, Baxter
  permutations and plane bipolar orientations, {\em S{\'e}minaire Lotharingien
  de Combinatoire} {\bf 61} (2010), B61Ah.

\bibitem{bose1998pattern}
Prosenjit Bose, Jonathan~F Buss, and Anna Lubiw, Pattern matching for
  permutations, {\em Information Processing Letters} {\bf 65}(5) (1998),
  277--283.

\bibitem{bouvel2018semi}
Mathilde Bouvel, Veronica Guerrini, Andrew Rechnitzer, and Simone Rinaldi,
  Semi-baxter and strong-baxter: two relatives of the baxter sequence, {\em
  SIAM Journal on Discrete Mathematics} {\bf 32}(4) (2018), 2795--2819.

\bibitem{brightwell1992random}
Graham Brightwell, Random k-dimensional orders: Width and number of linear
  extensions, {\em Order} {\bf 9}(4) (1992), 333--342.

\bibitem{burrill2016tableau}
Sophie Burrill, Julien Courtiel, Eric Fusy, Stephen Melczer, and Marni Mishna,
  Tableau sequences, open diagrams, and baxter families, {\em European Journal
  of Combinatorics} {\bf 58} (2016), 144--165.

\bibitem{cardinal_open}
Jean Cardinal, Order and geometry workshop 2016 - problem booklet,
  \url{http://orderandgeometry2016.tcs.uj.edu.pl/docs/OG2016-Problem-Booklet.pdf},
  2016.

\bibitem{silveira2018note}
Jean Cardinal, Vera Sacristán, and Rodrigo~I. Silveira, {A Note on Flips in
  Diagonal Rectangulations}, {\em {Discrete Mathematics \& Theoretical Computer
  Science}} {\bf {vol. 20 no. 2}} (2018), 1--22.

\bibitem{dulucq1998baxter}
Serge Dulucq and Olivier Guibert, Baxter permutations, {\em Discrete
  Mathematics} {\bf 180}(1-3) (1998), 143--156.

\bibitem{earnest2014permutation}
Michael~J Earnest and Samuel~C Gutekunst, Permutation patterns in {L}atin
  squares, {\em Australian Journoal of Combinatorics} {\bf 59}(1) (2014),
  218--228.

\bibitem{eriksson2000combinatorial}
Kimmo Eriksson and Svante Linusson, A combinatorial theory of
  higher-dimensional permutation arrays, {\em Advances in Applied Mathematics}
  {\bf 25}(2) (2000), 194--211.

\bibitem{felsner2011bijections}
Stefan Felsner, {\'E}ric Fusy, Marc Noy, and David Orden, Bijections for
  {B}axter families and related objects, {\em Journal of Combinatorial Theory,
  Series A} {\bf 118}(3) (2011), 993--1020.

\bibitem{felsner2020plattenbauten}
Stefan Felsner, Kolja Knauer, and Torsten Ueckerdt, Plattenbauten: touching
  rectangles in space, In {\em International Workshop on Graph-Theoretic
  Concepts in Computer Science}, pp.  161--173. Springer, 2020.

\bibitem{fusy_meanders}
Eric Fusy, Baxter permutations and meanders,
  \url{https://igm.univ-mlv.fr/~fusy/Talks/baxter_meanders.pdf}, 2012.

\bibitem{fusy2012bijective}
Eric Fusy, Bijective counting of involutive baxter permutations, {\em
  Fundamenta Informaticae} {\bf 117}(1-4) (2012), 179--188.

\bibitem{guibert2000doubly}
Olivier Guibert and Svante Linusson, Doubly alternating baxter permutations are
  catalan, {\em Discrete Mathematics} {\bf 217}(1-3) (2000), 157--166.

\bibitem{gunby2019asymptotics}
Benjamin Gunby and D{\"o}m{\"o}t{\"o}r P{\'a}lv{\"o}lgyi, Asymptotics of
  pattern avoidance in the {K}lazar set partition and permutation-tuple
  settings, {\em European Journal of Combinatorics} {\bf 82} (2019), 102992.

\bibitem{knuth1973art}
Donald~E Knuth, {\em The art of computer programming, volume 3: Searching and
  sorting}, Addison-Westley Publishing Company: Reading, MA, 1973.

\bibitem{mansour2020enumeration}
Toufik Mansour, Enumeration and {W}ilf-classification of permutations avoiding
  four patterns of length 4, {\em Discrete Math. Lett} {\bf 3} (2020), 67--94.

\bibitem{rotem1981stack}
Doron Rotem, Stack sortable permutations, {\em Discrete Mathematics} {\bf
  33}(2) (1981), 185--196.

\bibitem{salo2020cutting}
Ville Salo, Cutting corners, 2020.

\bibitem{simion1985restricted}
Rodica Simion and Frank~W Schmidt, Restricted permutations, {\em European
  Journal of Combinatorics} {\bf 6}(4) (1985), 383--406.

\bibitem{oeis}
Neil J.~A. Sloane and The OEIS~Foundation Inc., The on-line encyclopedia of
  integer sequences, 2020.

\bibitem{stembridge1996fully}
John~R Stembridge, On the fully commutative elements of {C}oxeter groups, {\em
  Journal of Algebraic Combinatorics} {\bf 5}(4) (1996), 353--385.

\bibitem{west2006enumeration}
Julian West and Olivier Guibert, Enumeration of reading’s twisted baxter
  permutations, In {\em The Fourth Annual International Conference on
  Permutation Patterns. Reykjavik University, Reykjavik}, 2006.

\end{thebibliography}
\bibliographystyle{jis}
\newpage
\appendix
\section{All symmetries of Baxter patterns}\label{sec:sym_bx}

$\Sym(\vinpatb{2413}{2}{2}) =$
$\vinpatb{2413}{2}{2},$
$\vinpatb{3142}{2}{2}.$

$\sbaxpb=\vinpatd{312}{213}{1}{.2}{.}$,
$\vinpatd{312}{231}{1}{.2}{.}$,
$\vinpatd{132}{213}{1}{.1}{.}$,
$\vinpatd{132}{231}{1}{.1}{.}$,
$\vinpatd{213}{312}{2}{.2}{.}$,
$\vinpatd{213}{132}{2}{.2}{.}$,
$\vinpatd{231}{312}{2}{.1}{.}$,
$\vinpatd{231}{132}{2}{.1}{.}$,
$\vinpatd{213}{312}{1}{.}{2}$,
$\vinpatd{213}{132}{1}{.}{1}$,
$\vinpatd{231}{312}{1}{.}{2}$,
$\vinpatd{231}{132}{1}{.}{1}$,
$\vinpatd{312}{213}{2}{.}{2}$,
$\vinpatd{312}{231}{2}{.}{1}$,
$\vinpatd{132}{213}{2}{.}{2}$,
$\vinpatd{132}{231}{2}{.}{1}$,
$\vinpatd{213}{132}{.}{1}{2}$,
$\vinpatd{213}{312}{.}{1}{1}$,
$\vinpatd{231}{132}{.}{2}{2}$,
$\vinpatd{231}{312}{.}{2}{1}$,
$\vinpatd{312}{231}{.}{1}{2}$,
$\vinpatd{312}{213}{.}{1}{1}$,
$\vinpatd{132}{231}{.}{2}{2}$,
$\vinpatd{132}{213}{.}{2}{1}$.

$\Sym(\vinpatd{3412}{1432}{2}{2}{.}) =$
$\vinpatd{2341}{4123}{.}{2}{2},$
$\vinpatd{2143}{3214}{2}{2}{.},$
$\vinpatd{4123}{3214}{.}{2}{2},$
$\vinpatd{3412}{3214}{2}{2}{.},$
$\vinpatd{3214}{4123}{.}{2}{2},$
$\vinpatd{2341}{1432}{.}{2}{2},$
$\vinpatd{1432}{3214}{.}{2}{2},$
$\vinpatd{2143}{1432}{2}{2}{.},$
$\vinpatd{3412}{1432}{2}{2}{.},$
$\vinpatd{2143}{4123}{2}{2}{.},$
$\vinpatd{1432}{2143}{2}{.}{2},$
$\vinpatd{4123}{2341}{.}{2}{2},$
$\vinpatd{3214}{1432}{.}{2}{2},$
$\vinpatd{3412}{4123}{2}{2}{.},$
$\vinpatd{3412}{2341}{2}{2}{.},$
$\vinpatd{1432}{3412}{2}{.}{2},$
$\vinpatd{2143}{2341}{2}{2}{.},$
$\vinpatd{2341}{3412}{2}{.}{2},$
$\vinpatd{4123}{2143}{2}{.}{2},$
$\vinpatd{4123}{3412}{2}{.}{2},$
$\vinpatd{3214}{3412}{2}{.}{2},$
$\vinpatd{1432}{2341}{.}{2}{2},$
$\vinpatd{3214}{2143}{2}{.}{2},$
$\vinpatd{2341}{2143}{2}{.}{2}.$

$\Sym(\vinpatd{2143}{1423}{2}{2}{.}) =$
$\vinpatd{3241}{2143}{2}{.}{2}, $
$\vinpatd{3412}{2314}{2}{2}{.}, $
$\vinpatd{1423}{3412}{2}{.}{2}, $
$\vinpatd{2314}{2143}{2}{.}{2}, $
$\vinpatd{1342}{3124}{.}{2}{2}, $
$\vinpatd{3124}{1342}{.}{2}{2}, $
$\vinpatd{1342}{2431}{.}{2}{2}, $
$\vinpatd{3241}{3412}{2}{.}{2}, $
$\vinpatd{4132}{3412}{2}{.}{2}, $
$\vinpatd{2431}{4213}{.}{2}{2}, $
$\vinpatd{2143}{3241}{2}{2}{.}, $
$\vinpatd{4213}{2431}{.}{2}{2}, $
$\vinpatd{3412}{3241}{2}{2}{.}, $
$\vinpatd{3412}{1423}{2}{2}{.}, $
$\vinpatd{4213}{3124}{.}{2}{2}, $
$\vinpatd{2143}{4132}{2}{2}{.}, $
$\vinpatd{3124}{4213}{.}{2}{2}, $
$\vinpatd{2431}{1342}{.}{2}{2}, $
$\vinpatd{2314}{3412}{2}{.}{2}, $
$\vinpatd{2143}{1423}{2}{2}{.}, $
$\vinpatd{1423}{2143}{2}{.}{2}, $
$\vinpatd{4132}{2143}{2}{.}{2}, $
$\vinpatd{2143}{2314}{2}{2}{.}, $
$\vinpatd{3412}{4132}{2}{2}{.}.$

\section{Other patterns}

We give here the beginning of sequences of permutations avoiding some larger
patterns or combination of patterns.
\begin{longtable}{|c|c|c|c|}
	\hline
	Patterns & \#TWE & Sequence & Comment \\
	\hline
	\hline
	1234 & 1 & 1, 4, 36, 506, 9032, 181582, 3836372, \ensuremath{\cdots} & \ensuremath{new} \\
\hline
1243 & 2 & 1, 4, 36, 507, 9089, 185253, 4017231, \ensuremath{\cdots} & \ensuremath{new} \\
\hline
1324 & 1 & 1, 4, 36, 507, 9087, 185455, 4053668, \ensuremath{\cdots} & \ensuremath{new} \\
\hline
1342 & 4 & 1, 4, 36, 507, 9102, 185920, 4059355, \ensuremath{\cdots} & \ensuremath{new} \\
\hline
1432 & 2 & 1, 4, 36, 507, 9119, 188501, 4230523, \ensuremath{\cdots} & \ensuremath{new} \\
\hline
2143 & 1 & 1, 4, 36, 507, 9121, 187799, 4163067, \ensuremath{\cdots} & \ensuremath{new} \\
\hline
2341 & 2 & 1, 4, 36, 507, 9105, 187502, 4191192, \ensuremath{\cdots} & \ensuremath{new} \\
\hline
2413 & 2 & 1, 4, 36, 507, 9141, 189810, 4291658, \ensuremath{\cdots} & \ensuremath{new} \\
\hline
2431 & 4 & 1, 4, 36, 507, 9124, 188197, 4197349, \ensuremath{\cdots} & \ensuremath{new} \\
\hline
3412 & 1 & 1, 4, 36, 507, 9135, 190457, 4368455, \ensuremath{\cdots} & \ensuremath{new} \\
\hline
3421 & 2 & 1, 4, 36, 507, 9133, 190307, 4355801, \ensuremath{\cdots} & \ensuremath{new} \\
\hline
4231 & 1 & 1, 4, 36, 507, 9119, 189363, 4318292, \ensuremath{\cdots} & \ensuremath{new} \\
\hline
4321 & 1 & 1, 4, 36, 507, 9147, 192181, 4482267, \ensuremath{\cdots} & \ensuremath{new} \\
\hline

	\caption{Patterns of size 4 and dimension 2.}
\end{longtable}

	\begin{longtable}{|c|c|c|c|}
		\hline
		Patterns & \#TWE & Sequence & Comment \\
		\hline
		\hline
		1234, 1243 & 2 & 1, 4, 36, 440, 5880, 75968, \ensuremath{\cdots} & \ensuremath{new} \\
\hline
1234, 1324 & 1 & 1, 4, 36, 440, 5872, 77616, \ensuremath{\cdots} & \ensuremath{new} \\
\hline
1234, 1342 & 4 & 1, 4, 36, 441, 5692, 68500, \ensuremath{\cdots} & \ensuremath{new} \\
\hline
1234, 1432 & 2 & 1, 4, 36, 440, 5056, 46446, \ensuremath{\cdots} & \ensuremath{new} \\
\hline
1234, 2143 & 1 & 1, 4, 36, 440, 5064, 45030, \ensuremath{\cdots} & \ensuremath{new} \\
\hline
1234, 2341 & 2 & 1, 4, 36, 441, 5730, 68040, \ensuremath{\cdots} & \ensuremath{new} \\
\hline
1234, 2413 & 2 & 1, 4, 36, 441, 5173, 49501, \ensuremath{\cdots} & \ensuremath{new} \\
\hline
1234, 2431 & 4 & 1, 4, 36, 441, 5180, 46360, \ensuremath{\cdots} & \ensuremath{new} \\
\hline
1234, 3412 & 1 & 1, 4, 36, 440, 5096, 44026, \ensuremath{\cdots} & \ensuremath{new} \\
\hline
1234, 3421 & 2 & 1, 4, 36, 441, 5205, 42991, \ensuremath{\cdots} & \ensuremath{new} \\
\hline
1234, 4231 & 1 & 1, 4, 36, 440, 5068, 43906, \ensuremath{\cdots} & \ensuremath{new} \\
\hline
1234, 4321 & 1 & 1, 4, 36, 440, 5168, 34784, \ensuremath{\cdots} & \ensuremath{new} \\
\hline
1243, 1324 & 2 & 1, 4, 36, 444, 6002, 79964, \ensuremath{\cdots} & \ensuremath{new} \\
\hline
1243, 1342 & 4 & 1, 4, 36, 444, 6015, 81001, \ensuremath{\cdots} & \ensuremath{new} \\
\hline
1243, 1432 & 2 & 1, 4, 36, 444, 5817, 73686, \ensuremath{\cdots} & \ensuremath{new} \\
\hline
1243, 2134 & 1 & 1, 4, 36, 444, 5353, 53256, \ensuremath{\cdots} & \ensuremath{new} \\
\hline
1243, 2143 & 2 & 1, 4, 36, 444, 6060, 82396, \ensuremath{\cdots} & \ensuremath{new} \\
\hline
1243, 2314 & 4 & 1, 4, 36, 444, 5647, 65690, \ensuremath{\cdots} & \ensuremath{new} \\
\hline
1243, 2341 & 4 & 1, 4, 36, 444, 5649, 65566, \ensuremath{\cdots} & \ensuremath{new} \\
\hline
1243, 2413 & 4 & 1, 4, 36, 444, 5700, 69626, \ensuremath{\cdots} & \ensuremath{new} \\
\hline
1243, 2431 & 4 & 1, 4, 36, 444, 5679, 66392, \ensuremath{\cdots} & \ensuremath{new} \\
\hline
1243, 3214 & 2 & 1, 4, 36, 444, 5278, 51226, \ensuremath{\cdots} & \ensuremath{new} \\
\hline
1243, 3241 & 4 & 1, 4, 36, 444, 5339, 54622, \ensuremath{\cdots} & \ensuremath{new} \\
\hline
1243, 3412 & 2 & 1, 4, 36, 444, 5336, 54613, \ensuremath{\cdots} & \ensuremath{new} \\
\hline
1243, 3421 & 4 & 1, 4, 36, 444, 5336, 51612, \ensuremath{\cdots} & \ensuremath{new} \\
\hline
1243, 4231 & 2 & 1, 4, 36, 444, 5296, 52363, \ensuremath{\cdots} & \ensuremath{new} \\
\hline
1243, 4321 & 2 & 1, 4, 36, 444, 5324, 47835, \ensuremath{\cdots} & \ensuremath{new} \\
\hline
1324, 1342 & 4 & 1, 4, 36, 444, 6036, 82584, \ensuremath{\cdots} & \ensuremath{new} \\
\hline
1324, 1432 & 2 & 1, 4, 36, 444, 5827, 73608, \ensuremath{\cdots} & \ensuremath{new} \\
\hline
1324, 2143 & 1 & 1, 4, 36, 444, 5650, 65194, \ensuremath{\cdots} & \ensuremath{new} \\
\hline
1324, 2341 & 2 & 1, 4, 36, 444, 5468, 59406, \ensuremath{\cdots} & \ensuremath{new} \\
\hline
1324, 2413 & 2 & 1, 4, 36, 444, 5726, 70540, \ensuremath{\cdots} & \ensuremath{new} \\
\hline
1324, 2431 & 4 & 1, 4, 36, 444, 5710, 68014, \ensuremath{\cdots} & \ensuremath{new} \\
\hline
1324, 3412 & 1 & 1, 4, 36, 444, 5304, 52359, \ensuremath{\cdots} & \ensuremath{new} \\
\hline
1324, 3421 & 2 & 1, 4, 36, 444, 5317, 53022, \ensuremath{\cdots} & \ensuremath{new} \\
\hline
1324, 4231 & 1 & 1, 4, 36, 444, 5276, 52016, \ensuremath{\cdots} & \ensuremath{new} \\
\hline
1324, 4321 & 1 & 1, 4, 36, 444, 5304, 50792, \ensuremath{\cdots} & \ensuremath{new} \\
\hline
1342, 1423 & 2 & 1, 4, 36, 442, 5978, 82076, \ensuremath{\cdots} & \ensuremath{new} \\
\hline
1342, 1432 & 4 & 1, 4, 36, 444, 6056, 84402, \ensuremath{\cdots} & \ensuremath{new} \\
\hline
1342, 2143 & 4 & 1, 4, 36, 444, 5692, 68333, \ensuremath{\cdots} & \ensuremath{new} \\
\hline
1342, 2314 & 2 & 1, 4, 36, 444, 5710, 69187, \ensuremath{\cdots} & \ensuremath{new} \\
\hline
1342, 2341 & 4 & 1, 4, 36, 444, 6080, 84954, \ensuremath{\cdots} & \ensuremath{new} \\
\hline
1342, 2413 & 4 & 1, 4, 36, 444, 5952, 80102, \ensuremath{\cdots} & \ensuremath{new} \\
\hline
1342, 2431 & 4 & 1, 4, 36, 444, 5726, 70904, \ensuremath{\cdots} & \ensuremath{new} \\
\hline
1342, 3124 & 2 & 1, 4, 36, 444, 5507, 62078, \ensuremath{\cdots} & \ensuremath{new} \\
\hline
1342, 3142 & 4 & 1, 4, 36, 444, 6148, 88944, \ensuremath{\cdots} & \ensuremath{new} \\
\hline
1342, 3214 & 4 & 1, 4, 36, 444, 5334, 54125, \ensuremath{\cdots} & \ensuremath{new} \\
\hline
1342, 3241 & 4 & 1, 4, 36, 444, 5733, 70753, \ensuremath{\cdots} & \ensuremath{new} \\
\hline
1342, 3412 & 4 & 1, 4, 36, 444, 5738, 71301, \ensuremath{\cdots} & \ensuremath{new} \\
\hline
1342, 3421 & 4 & 1, 4, 36, 444, 5715, 68527, \ensuremath{\cdots} & \ensuremath{new} \\
\hline
1342, 4123 & 4 & 1, 4, 36, 444, 5483, 60355, \ensuremath{\cdots} & \ensuremath{new} \\
\hline
1342, 4132 & 4 & 1, 4, 36, 444, 5734, 70864, \ensuremath{\cdots} & \ensuremath{new} \\
\hline
1342, 4213 & 4 & 1, 4, 36, 444, 5364, 56948, \ensuremath{\cdots} & \ensuremath{new} \\
\hline
1342, 4231 & 4 & 1, 4, 36, 444, 5706, 68457, \ensuremath{\cdots} & \ensuremath{new} \\
\hline
1342, 4312 & 4 & 1, 4, 36, 444, 5356, 56450, \ensuremath{\cdots} & \ensuremath{new} \\
\hline
1342, 4321 & 4 & 1, 4, 36, 444, 5324, 51799, \ensuremath{\cdots} & \ensuremath{new} \\
\hline
1432, 2143 & 2 & 1, 4, 36, 444, 5931, 77775, \ensuremath{\cdots} & \ensuremath{new} \\
\hline
1432, 2341 & 4 & 1, 4, 36, 444, 5348, 57776, \ensuremath{\cdots} & \ensuremath{new} \\
\hline
1432, 2413 & 4 & 1, 4, 36, 444, 5766, 73833, \ensuremath{\cdots} & \ensuremath{new} \\
\hline
1432, 2431 & 4 & 1, 4, 36, 444, 6126, 87630, \ensuremath{\cdots} & \ensuremath{new} \\
\hline
1432, 3214 & 1 & 1, 4, 36, 444, 5587, 63160, \ensuremath{\cdots} & \ensuremath{new} \\
\hline
1432, 3241 & 4 & 1, 4, 36, 444, 5536, 63590, \ensuremath{\cdots} & \ensuremath{new} \\
\hline
1432, 3412 & 2 & 1, 4, 36, 444, 5444, 63144, \ensuremath{\cdots} & \ensuremath{new} \\
\hline
1432, 3421 & 4 & 1, 4, 36, 444, 5761, 72105, \ensuremath{\cdots} & \ensuremath{new} \\
\hline
1432, 4231 & 2 & 1, 4, 36, 444, 5485, 62074, \ensuremath{\cdots} & \ensuremath{new} \\
\hline
1432, 4321 & 2 & 1, 4, 36, 444, 5981, 79272, \ensuremath{\cdots} & \ensuremath{new} \\
\hline
2143, 2341 & 2 & 1, 4, 36, 444, 5349, 56637, \ensuremath{\cdots} & \ensuremath{new} \\
\hline
2143, 2413 & 2 & 1, 4, 36, 444, 6146, 88824, \ensuremath{\cdots} & \ensuremath{new} \\
\hline
2143, 2431 & 4 & 1, 4, 36, 444, 5730, 70097, \ensuremath{\cdots} & \ensuremath{new} \\
\hline
2143, 3412 & 1 & 1, 4, 36, 444, 5476, 62504, \ensuremath{\cdots} & \ensuremath{new} \\
\hline
2143, 3421 & 2 & 1, 4, 36, 443, 5357, 56583, \ensuremath{\cdots} & \ensuremath{new} \\
\hline
2143, 4231 & 1 & 1, 4, 36, 444, 5322, 53529, \ensuremath{\cdots} & \ensuremath{new} \\
\hline
2143, 4321 & 1 & 1, 4, 36, 444, 5464, 58437, \ensuremath{\cdots} & \ensuremath{new} \\
\hline
2341, 2413 & 4 & 1, 4, 36, 444, 5731, 72541, \ensuremath{\cdots} & \ensuremath{new} \\
\hline
2341, 2431 & 4 & 1, 4, 36, 444, 6122, 87944, \ensuremath{\cdots} & \ensuremath{new} \\
\hline
2341, 3412 & 2 & 1, 4, 36, 443, 5864, 77512, \ensuremath{\cdots} & \ensuremath{new} \\
\hline
2341, 3421 & 2 & 1, 4, 36, 444, 5922, 80471, \ensuremath{\cdots} & \ensuremath{new} \\
\hline
2341, 4123 & 1 & 1, 4, 36, 444, 5441, 56318, \ensuremath{\cdots} & \ensuremath{new} \\
\hline
2341, 4132 & 4 & 1, 4, 36, 444, 5329, 54619, \ensuremath{\cdots} & \ensuremath{new} \\
\hline
2341, 4231 & 2 & 1, 4, 36, 444, 5894, 78113, \ensuremath{\cdots} & \ensuremath{new} \\
\hline
2341, 4312 & 2 & 1, 4, 36, 444, 5342, 56655, \ensuremath{\cdots} & \ensuremath{new} \\
\hline
2341, 4321 & 2 & 1, 4, 36, 444, 5371, 60374, \ensuremath{\cdots} & \ensuremath{new} \\
\hline
2413, 2431 & 4 & 1, 4, 36, 444, 6164, 89724, \ensuremath{\cdots} & \ensuremath{new} \\
\hline
2413, 3142 & 1 & 1, 4, 36, 444, 6252, 94588, \ensuremath{\cdots} & \ensuremath{new} \\
\hline
2413, 3241 & 4 & 1, 4, 36, 444, 5962, 80566, \ensuremath{\cdots} & \ensuremath{new} \\
\hline
2413, 3412 & 2 & 1, 4, 36, 444, 6162, 90477, \ensuremath{\cdots} & \ensuremath{new} \\
\hline
2413, 3421 & 4 & 1, 4, 36, 444, 5746, 72759, \ensuremath{\cdots} & \ensuremath{new} \\
\hline
2413, 4231 & 2 & 1, 4, 36, 444, 5760, 72775, \ensuremath{\cdots} & \ensuremath{new} \\
\hline
2413, 4321 & 2 & 1, 4, 36, 443, 5359, 58000, \ensuremath{\cdots} & \ensuremath{new} \\
\hline
2431, 3241 & 2 & 1, 4, 36, 444, 6137, 88439, \ensuremath{\cdots} & \ensuremath{new} \\
\hline
2431, 3412 & 4 & 1, 4, 36, 444, 5758, 73920, \ensuremath{\cdots} & \ensuremath{new} \\
\hline
2431, 3421 & 4 & 1, 4, 36, 444, 6149, 89342, \ensuremath{\cdots} & \ensuremath{new} \\
\hline
2431, 4132 & 2 & 1, 4, 36, 442, 5662, 70024, \ensuremath{\cdots} & \ensuremath{new} \\
\hline
2431, 4213 & 2 & 1, 4, 36, 444, 5565, 65925, \ensuremath{\cdots} & \ensuremath{new} \\
\hline
2431, 4231 & 4 & 1, 4, 36, 444, 6134, 88594, \ensuremath{\cdots} & \ensuremath{new} \\
\hline
2431, 4312 & 4 & 1, 4, 36, 444, 5754, 73295, \ensuremath{\cdots} & \ensuremath{new} \\
\hline
2431, 4321 & 4 & 1, 4, 36, 444, 5978, 82140, \ensuremath{\cdots} & \ensuremath{new} \\
\hline
3412, 3421 & 2 & 1, 4, 36, 444, 6196, 91640, \ensuremath{\cdots} & \ensuremath{new} \\
\hline
3412, 4231 & 1 & 1, 4, 36, 444, 5726, 72248, \ensuremath{\cdots} & \ensuremath{new} \\
\hline
3412, 4321 & 1 & 1, 4, 36, 444, 5496, 66138, \ensuremath{\cdots} & \ensuremath{new} \\
\hline
3421, 4231 & 2 & 1, 4, 36, 444, 6152, 90102, \ensuremath{\cdots} & \ensuremath{new} \\
\hline
3421, 4312 & 1 & 1, 4, 36, 444, 5655, 70866, \ensuremath{\cdots} & \ensuremath{new} \\
\hline
3421, 4321 & 2 & 1, 4, 36, 444, 6228, 93468, \ensuremath{\cdots} & \ensuremath{new} \\
\hline
4231, 4321 & 1 & 1, 4, 36, 444, 6176, 92820, \ensuremath{\cdots} & \ensuremath{new} \\
\hline

		\caption{Pairs of patterns of size 4 and dimension 2.}
	\end{longtable}

\begin{longtable}{|c|c|c|c|}
	\hline
	Patterns & \#TWE & Sequence & Comment \\
	\hline
	\hline
	123, \perm{123}{123} & 1 & 1, 4, 20, 100, 410, 1224, 2232, \ensuremath{\cdots} & 123 \\
\hline
123, \perm{123}{132} & 6 & 1, 4, 20, 100, 410, 1224, 2232, \ensuremath{\cdots} & 123 \\
\hline
123, \perm{123}{231} & 6 & 1, 4, 20, 100, 410, 1224, 2232, \ensuremath{\cdots} & 123 \\
\hline
123, \perm{123}{321} & 3 & 1, 4, 20, 100, 410, 1224, 2232, \ensuremath{\cdots} & 123 \\
\hline
123, \perm{132}{213} & 6 & 1, 4, 19, 91, 358, 1005, 1601, \ensuremath{\cdots} & \ensuremath{new} \\
\hline
123, \perm{132}{312} & 12 & 1, 4, 19, 79, 231, 407, 354, \ensuremath{\cdots} & \ensuremath{new} \\
\hline
123, \perm{231}{312} & 2 & 1, 4, 19, 83, 262, 514, 527, \ensuremath{\cdots} & \ensuremath{new} \\
\hline
132, \perm{123}{123} & 2 & 1, 4, 20, 100, 490, 2366, 11334, \ensuremath{\cdots} & \ensuremath{new} \\
\hline
132, \perm{123}{132} & 6 & 1, 4, 21, 116, 646, 3596, 19981, \ensuremath{\cdots} & 132 \\
\hline
132, \perm{123}{213} & 6 & 1, 4, 20, 102, 518, 2618, 13194, \ensuremath{\cdots} & \ensuremath{new} \\
\hline
132, \perm{123}{231} & 6 & 1, 4, 20, 100, 486, 2302, 10690, \ensuremath{\cdots} & \ensuremath{new} \\
\hline
132, \perm{123}{312} & 6 & 1, 4, 20, 104, 544, 2846, 14880, \ensuremath{\cdots} & \ensuremath{new} \\
\hline
132, \perm{123}{321} & 6 & 1, 4, 20, 99, 477, 2252, 10480, \ensuremath{\cdots} & \ensuremath{new} \\
\hline
132, \perm{132}{213} & 12 & 1, 4, 21, 116, 646, 3596, 19981, \ensuremath{\cdots} & 132 \\
\hline
132, \perm{132}{312} & 12 & 1, 4, 21, 116, 646, 3596, 19981, \ensuremath{\cdots} & 132 \\
\hline
132, \perm{213}{231} & 12 & 1, 4, 20, 100, 488, 2335, 11016, \ensuremath{\cdots} & \ensuremath{new} \\
\hline
132, \perm{231}{312} & 4 & 1, 4, 20, 105, 559, 2990, 16021, \ensuremath{\cdots} & \ensuremath{new} \\
\hline
231, \perm{123}{123} & 2 & 1, 4, 20, 97, 431, 1758, 6669, \ensuremath{\cdots} & \ensuremath{new} \\
\hline
231, \perm{123}{132} & 4 & 1, 4, 20, 104, 544, 2855, 15056, \ensuremath{\cdots} & \ensuremath{new} \\
\hline
231, \perm{123}{213} & 4 & 1, 4, 20, 106, 573, 3127, 17173, \ensuremath{\cdots} & \ensuremath{new} \\
\hline
231, \perm{123}{231} & 4 & 1, 4, 21, 123, 767, 4994, 33584, \ensuremath{\cdots} & 231 \\
\hline
231, \perm{123}{312} & 4 & 1, 4, 20, 105, 564, 3094, 17329, \ensuremath{\cdots} & \ensuremath{new} \\
\hline
231, \perm{123}{321} & 4 & 1, 4, 20, 106, 581, 3273, 18851, \ensuremath{\cdots} & \ensuremath{new} \\
\hline
231, \perm{132}{123} & 4 & 1, 4, 20, 105, 564, 3092, 17289, \ensuremath{\cdots} & \ensuremath{new} \\
\hline
231, \perm{132}{213} & 4 & 1, 4, 21, 123, 767, 4994, 33584, \ensuremath{\cdots} & 231 \\
\hline
231, \perm{132}{231} & 2 & 1, 4, 21, 123, 767, 4994, 33584, \ensuremath{\cdots} & 231 \\
\hline
231, \perm{132}{312} & 4 & 1, 4, 20, 108, 611, 3575, 21455, \ensuremath{\cdots} & \ensuremath{new} \\
\hline
231, \perm{132}{321} & 4 & 1, 4, 20, 108, 607, 3504, 20638, \ensuremath{\cdots} & \ensuremath{new} \\
\hline
231, \perm{213}{132} & 4 & 1, 4, 20, 109, 629, 3793, 23669, \ensuremath{\cdots} & \ensuremath{new} \\
\hline
231, \perm{213}{231} & 4 & 1, 4, 21, 123, 767, 4994, 33584, \ensuremath{\cdots} & 231 \\
\hline
231, \perm{213}{312} & 2 & 1, 4, 20, 111, 654, 4013, 25380, \ensuremath{\cdots} & \ensuremath{new} \\
\hline
231, \perm{213}{321} & 4 & 1, 4, 21, 123, 767, 4994, 33584, \ensuremath{\cdots} & 231 \\
\hline
231, \perm{231}{123} & 4 & 1, 4, 21, 123, 767, 4994, 33584, \ensuremath{\cdots} & 231 \\
\hline
231, \perm{231}{213} & 4 & 1, 4, 21, 123, 767, 4994, 33584, \ensuremath{\cdots} & 231 \\
\hline
231, \perm{231}{312} & 2 & 1, 4, 21, 123, 767, 4994, 33584, \ensuremath{\cdots} & 231 \\
\hline
231, \perm{312}{132} & 4 & 1, 4, 20, 111, 659, 4102, 26435, \ensuremath{\cdots} & \ensuremath{new} \\
\hline
231, \perm{312}{231} & 2 & 1, 4, 21, 123, 767, 4994, 33584, \ensuremath{\cdots} & 231 \\
\hline
231, \perm{321}{123} & 2 & 1, 4, 20, 112, 673, 4243, 27696, \ensuremath{\cdots} & \ensuremath{new} \\
\hline
321, \perm{123}{123} & 1 & 1, 4, 20, 76, 108, 52, 0, \ensuremath{\cdots} &  \\
\hline
321, \perm{123}{132} & 6 & 1, 4, 20, 103, 527, 2714, 14274, \ensuremath{\cdots} & \ensuremath{new} \\
\hline
321, \perm{123}{231} & 6 & 1, 4, 20, 110, 644, 3934, 24770, \ensuremath{\cdots} & \ensuremath{new} \\
\hline
321, \perm{123}{321} & 3 & 1, 4, 21, 128, 850, 5956, 43235, \ensuremath{\cdots} & 321 \\
\hline
321, \perm{132}{213} & 6 & 1, 4, 20, 113, 687, 4389, 29046, \ensuremath{\cdots} & \ensuremath{new} \\
\hline
321, \perm{132}{312} & 12 & 1, 4, 21, 128, 850, 5956, 43235, \ensuremath{\cdots} & 321 \\
\hline
321, \perm{231}{312} & 2 & 1, 4, 20, 117, 745, 5006, 34873, \ensuremath{\cdots} & \ensuremath{new} \\
\hline

	\caption{Pairs of patterns of size 3 respectively of dimension 2 and 3. }
\end{longtable}

\begin{longtable}{|c|c|c|c|}
	\hline
	Patterns & \#TWE & Sequence & Comment \\
	\hline
	\hline
	\perm{123}{123}, \perm{123}{132} & 24 & 1, 4, 34, 480, 9916, 277730, 10023010, \ensuremath{\cdots} & \ensuremath{new} \\
\hline
\perm{123}{123}, \perm{123}{231} & 24 & 1, 4, 34, 477, 9681, 262606, 9038034, \ensuremath{\cdots} & \ensuremath{new} \\
\hline
\perm{123}{123}, \perm{123}{321} & 6 & 1, 4, 34, 472, 9324, 241616, 7793548, \ensuremath{\cdots} & \ensuremath{new} \\
\hline
\perm{123}{123}, \perm{132}{213} & 24 & 1, 4, 34, 476, 9618, 259274, 8857074, \ensuremath{\cdots} & \ensuremath{new} \\
\hline
\perm{123}{123}, \perm{132}{312} & 48 & 1, 4, 34, 472, 9321, 241306, 7769550, \ensuremath{\cdots} & \ensuremath{new} \\
\hline
\perm{123}{123}, \perm{231}{312} & 8 & 1, 4, 34, 472, 9286, 237532, 7466512, \ensuremath{\cdots} & \ensuremath{new} \\
\hline
\perm{123}{132}, \perm{123}{213} & 12 & 1, 4, 34, 478, 9758, 267578, 9366032, \ensuremath{\cdots} & \ensuremath{new} \\
\hline
\perm{123}{132}, \perm{123}{231} & 12 & 1, 4, 34, 480, 9916, 277792, 10032960, \ensuremath{\cdots} & \ensuremath{new} \\
\hline
\perm{123}{132}, \perm{123}{312} & 12 & 1, 4, 34, 476, 9622, 259720, 8895656, \ensuremath{\cdots} & \ensuremath{new} \\
\hline
\perm{123}{132}, \perm{132}{123} & 24 & 1, 4, 34, 480, 9912, 277304, 9987248, \ensuremath{\cdots} & \ensuremath{new} \\
\hline
\perm{123}{132}, \perm{132}{213} & 48 & 1, 4, 34, 476, 9617, 259152, 8846076, \ensuremath{\cdots} & \ensuremath{new} \\
\hline
\perm{123}{132}, \perm{132}{312} & 48 & 1, 4, 34, 474, 9463, 249551, 8249751, \ensuremath{\cdots} & \ensuremath{new} \\
\hline
\perm{123}{132}, \perm{213}{123} & 24 & 1, 4, 34, 476, 9633, 260990, 9007402, \ensuremath{\cdots} & \ensuremath{new} \\
\hline
\perm{123}{132}, \perm{213}{132} & 48 & 1, 4, 34, 480, 9900, 275992, 9874628, \ensuremath{\cdots} & \ensuremath{new} \\
\hline
\perm{123}{132}, \perm{213}{231} & 48 & 1, 4, 34, 475, 9555, 255962, 8679070, \ensuremath{\cdots} & \ensuremath{new} \\
\hline
\perm{123}{132}, \perm{231}{132} & 48 & 1, 4, 34, 476, 9608, 258290, 8782799, \ensuremath{\cdots} & \ensuremath{new} \\
\hline
\perm{123}{132}, \perm{231}{213} & 24 & 1, 4, 34, 474, 9462, 249440, 8240370, \ensuremath{\cdots} & \ensuremath{new} \\
\hline
\perm{123}{132}, \perm{231}{312} & 48 & 1, 4, 34, 474, 9441, 247195, 8060190, \ensuremath{\cdots} & \ensuremath{new} \\
\hline
\perm{123}{132}, \perm{231}{321} & 24 & 1, 4, 34, 476, 9603, 257690, 8728931, \ensuremath{\cdots} & \ensuremath{new} \\
\hline
\perm{123}{132}, \perm{321}{132} & 24 & 1, 4, 34, 472, 9332, 242344, 7844248, \ensuremath{\cdots} & \ensuremath{new} \\
\hline
\perm{123}{132}, \perm{321}{213} & 24 & 1, 4, 34, 472, 9316, 240804, 7731538, \ensuremath{\cdots} & \ensuremath{new} \\
\hline
\perm{132}{213}, \perm{132}{231} & 12 & 1, 4, 34, 476, 9618, 259364, 8871444, \ensuremath{\cdots} & \ensuremath{new} \\
\hline
\perm{132}{213}, \perm{213}{132} & 4 & 1, 4, 34, 478, 9730, 264334, 9076864, \ensuremath{\cdots} & \ensuremath{new} \\
\hline
\perm{132}{213}, \perm{213}{312} & 12 & 1, 4, 34, 474, 9450, 248156, 8137074, \ensuremath{\cdots} & \ensuremath{new} \\
\hline

	\caption{Pairs of patterns of size 3 and of dimension 3.}
\end{longtable}

\begin{longtable}{|c|c|c|c|}
	\hline
	Patterns & \#TWE & Sequence & Comment \\
	\hline
	\hline
	\perm{1234}{1234} & 4 & 1, 4, 36, 575, 14291, 508161, 24385927, \ensuremath{\cdots} & \ensuremath{new} \\
\hline
\perm{1234}{1243} & 24 & 1, 4, 36, 575, 14291, 508155, 24384283, \ensuremath{\cdots} & \ensuremath{new} \\
\hline
\perm{1234}{1324} & 12 & 1, 4, 36, 575, 14291, 508149, 24382888, \ensuremath{\cdots} & \ensuremath{new} \\
\hline
\perm{1234}{1342} & 24 & 1, 4, 36, 575, 14291, 508144, 24381346, \ensuremath{\cdots} & \ensuremath{new} \\
\hline
\perm{1234}{1423} & 24 & 1, 4, 36, 575, 14291, 508144, 24381396, \ensuremath{\cdots} & \ensuremath{new} \\
\hline
\perm{1234}{1432} & 24 & 1, 4, 36, 575, 14291, 508155, 24384181, \ensuremath{\cdots} & \ensuremath{new} \\
\hline
\perm{1234}{2143} & 12 & 1, 4, 36, 575, 14291, 508153, 24383579, \ensuremath{\cdots} & \ensuremath{new} \\
\hline
\perm{1234}{2413} & 12 & 1, 4, 36, 575, 14291, 508132, 24378096, \ensuremath{\cdots} & \ensuremath{new} \\
\hline
\perm{1243}{1324} & 48 & 1, 4, 36, 575, 14291, 508135, 24379128, \ensuremath{\cdots} & \ensuremath{new} \\
\hline
\perm{1243}{1423} & 48 & 1, 4, 36, 575, 14291, 508144, 24381329, \ensuremath{\cdots} & \ensuremath{new} \\
\hline
\perm{1243}{2134} & 24 & 1, 4, 36, 575, 14291, 508151, 24383081, \ensuremath{\cdots} & \ensuremath{new} \\
\hline
\perm{1243}{2314} & 48 & 1, 4, 36, 575, 14291, 508142, 24380642, \ensuremath{\cdots} & \ensuremath{new} \\
\hline
\perm{1243}{2413} & 48 & 1, 4, 36, 575, 14291, 508129, 24377368, \ensuremath{\cdots} & \ensuremath{new} \\
\hline
\perm{1324}{1342} & 48 & 1, 4, 36, 575, 14291, 508142, 24380847, \ensuremath{\cdots} & \ensuremath{new} \\
\hline
\perm{1324}{2143} & 24 & 1, 4, 36, 575, 14291, 508131, 24377763, \ensuremath{\cdots} & \ensuremath{new} \\
\hline
\perm{1342}{1423} & 16 & 1, 4, 36, 575, 14291, 508131, 24378031, \ensuremath{\cdots} & \ensuremath{new} \\
\hline
\perm{1342}{2143} & 24 & 1, 4, 36, 575, 14291, 508132, 24378046, \ensuremath{\cdots} & \ensuremath{new} \\
\hline
\perm{1342}{2314} & 16 & 1, 4, 36, 575, 14291, 508128, 24377163, \ensuremath{\cdots} & \ensuremath{new} \\
\hline
\perm{1342}{2413} & 48 & 1, 4, 36, 575, 14291, 508128, 24377001, \ensuremath{\cdots} & \ensuremath{new} \\
\hline
\perm{1342}{2431} & 24 & 1, 4, 36, 575, 14291, 508139, 24379797, \ensuremath{\cdots} & \ensuremath{new} \\
\hline
\perm{1432}{2143} & 24 & 1, 4, 36, 575, 14291, 508143, 24380822, \ensuremath{\cdots} & \ensuremath{new} \\
\hline

	\caption{Patterns of size 4 and dimension 3.}
\end{longtable}

\end{document}